\documentclass[12pt]{article}
\usepackage{amsmath,amssymb,showlabels}

\newtheorem{Theorem}{Theorem}[part]
\newtheorem{Definition}{Definition}[part]
\newtheorem{Proposition}{Proposition}[part]
\newtheorem{Assumption}{Assumption}[part]

\newtheorem{Corollary}{Corollary}[part]

\def \H{I\!\!H}

\def \R{I\!\!R}

\def\Ac{{\cal A}}

\def\Fc{{\cal F}}

\def\T{{\cal T}}

\def\Dzw1#1{\frac{\partial^2 #1}{\partial z \partial w_1}}

\def\Dzb1#1{\frac{\partial^2 #1}{\partial z \partial b_1}}

\newcommand{\dproof}{\noindent {Proof.} \quad}
\newcommand{\fproof}{\hfill $\square$ \bigskip}

\newtheorem{definition}{Definition}[section]

\newtheorem{theorem}[definition]{Theorem}

\newtheorem{remark}[definition]{ \it Remark}

\newtheorem{lemma}[definition]{Lemma}

\def\RB{\mathbb{R}}

\def\EB{\mathbb{E}}

\def\BC{\mathcal{B}}
\def\FC{\mathcal{F}}
\def\PC{\mathcal{P}}

\def\R{{\bf R}}

\def\1B{\text{1\!\!I}}

\def\tN{\tilde{N}}

\def\RB{\mathbb{R}}

\def\EB{\mathbb{E}}

\def\FC{\mathcal{F}}
\def\PC{\mathcal{P}}

\def\R{{\bf R}}

\def\1B{\text{1\!\!I}}

\def\tN{\tilde{N}}

\def\RB{\mathbb{R}}

\def\EB{\mathbb{E}}

\begin{document}

\title{Reflected BSDEs and robust optimal stopping for dynamic risk measures with jumps} 
\author{Marie-Claire QUENEZ\thanks{LPMA, 
Universit\'e Paris 7 Denis Diderot, Boite courrier 7012, 75251 Paris cedex 05, France, and  INRIA Paris-Rocquencourt, Domaine de Voluceau, Rocquencourt, BP 105, Le Chesnay Cedex, 78153, France, email: {\tt quenez@math.univ-paris-diderot.fr}}
\and 
Agn\`es SULEM
\thanks{ INRIA Paris-Rocquencourt, Domaine de Voluceau, Rocquencourt, BP 105, Le Chesnay Cedex, 78153, France, and Universit\'e Paris-Est, email: {\tt agnes.sulem@inria.fr}}}

\date{\today}

\maketitle

\begin{abstract}
We  study  the optimal stopping   problem for  dynamic risk measures represented by Backward Stochastic Differential Equations (BSDEs)  with jumps and its relation with reflected BSDEs (RBSDEs). 
We first provide general existence, uniqueness and  comparison  theorems for RBSDEs with jumps in the case of a  RCLL adapted obstacle.
We  then show  that  the value function of the optimal  stopping problem  is characterized as the solution  of an RBSDE. 
 The existence  of an optimal stopping time is obtained   when the obstacle is left-upper semi-continuous along stopping times.
Finally, robust optimal stopping   problems  related to the case  with  model ambiguity are investigated.

%
%
\end{abstract}

\vspace{13mm}

\noindent{\bf Key words~:}  Backward stochastic differential equations, reflected backward stochastic equations,  jump processes,  optimal stopping,   risk-measures.

\vspace{13mm}

\noindent{\bf AMS 1991 subject classifications~:} 93E20, 60J60, 47N10.

\section{Introduction}\label{sec1}

In this paper we study   optimal stopping   problems for  dynamic risk measures  $\rho_t$ represented by Backward Stochastic Differential Equations (BSDEs)  with jumps. The properties of these  risk measures have been studied recently in \cite{QuenSul}. 
The optimal stopping problem can be formulated as follows: 
 given a dynamic financial position $\xi_t$, represented  by an RCLL adapted process, we want to determine a stopping time  $\tau$ which minimizes the risk of the position $\xi_\tau$,
and compute the corresponding value. 
%
To this purpose, we study the links between this optimal stopping problem and reflected BSDEs  (RBSDEs) with jumps.
RBSDEs  have been introduced by  N. El Karoui et al. (1997 ) (see \cite{ElKaroui97}) in the case of a Brownian filtration. The solutions of such equations 
are constrained to be greater than  given processes called obstacles.
We  provide here existence and uniqueness results for RBSDEs with jumps, as well as comparison and strict comparison  theorems,   when the obstacle is RCLL.  This  completes some results in Hamad\` ene, Ouknine and Issaky \cite{HO1, HO2, Essaky}.

We prove  that  the value function of our optimal stopping problem  is the solution  
of an RBSDE with obstacle given by the dynamic position $\xi_t$.  
We provide an optimality criterium, that is a characterization of optimal stopping times. In the case when the obstacle is left-upper semi-continuous along stopping times, 
we show the existence of an optimal stopping time. In the case of a general RCLL obstacle, 
we prove the existence of $\varepsilon$-stopping times. Related studies can be found in El Karoui and Quenez \cite{EQ96}, Bayraktar and coauthors in \cite{Bayraktar} and \cite{Bayraktar-2} in the Brownian case. 
  

 We then address the optimal stopping problem when there is ambiguity on the risk measure. To this purpose, we  study the following optimal control problem for RBSDEs: 
Let $\{f^{\alpha}, \alpha \in {\cal A}\}$ be a 
 family of Lipschitz drivers and let $\{Y^{\alpha}, \alpha \in {\cal A}\}$ be the solutions of the  RBSDEs associated 
 with drivers $\{f^{\alpha}\}$ and obstacle $\xi_t$. The problem is to minimize   $Y^{\alpha}$ over $\alpha$. Under appropriate hypotheses, the value function is characterized as the solution $Y$  of an RBSDE.  
We then focus on the robust optimal stopping problem for risk measures: 
 we consider the family of risk measures $\{\rho_t^{\alpha},\alpha \in {\cal A}\}$  induced by the 
BSDEs associated with  drivers 
$\{f^{\alpha}, \alpha \in {\cal A}\}$.  
In this ambiguity framework,  the risk measure is defined as the supremum over $\alpha$ of the risk measures $\rho^\alpha$.
Given the dynamic position $\xi_t$,  we want to determine a stopping time  $\tau^*$ which minimizes  over all stopping times $\tau$ the risk of the position $\xi_{\tau}.$ 
This leads to a mixed control/ optimal stopping game problem. 
We show 
that, under some hypothesis, the value function is equal to $Y $. 
We then study the existence of saddle points.

The paper is organized as follows. In Section~\ref{sec2}, we state the notation and give the formulation of our optimal stopping problem for risk measures.
In Section~\ref{rbsde}, we provide existence and uniqueness results for RBSDEs with jumps 
and RCLL obstacle. 
Relations between optimal stopping problems and   RBSDEs  are given in Section \ref{sec-charact}.
In Section~\ref{sec-charc}, we provide  comparison theorems for RBSDEs with jumps and optimization principles. 
The robust optimal stopping problem for risk measures when there is ambiguity on the risk measure is addressed in Section~\ref{mixed}. An application to a case of multiple priors is presented in Section~\ref{application}.

\section{Formulation of the problem} \label{sec2}


\paragraph{Notation.}
Let ${\cal P}$ be  the predictable $\sigma$-algebra
on $[0,T]  \times \Omega$. 

 For each $T>0$ and $p>1$, we use the following notation: 
\begin{itemize}
\item 
 $L^p({\cal F}_T)$  is the set of random variables $\xi$ which are  $\Fc 
_T$-measurable and $p$-integrable.

\item    $\H^{p,T}$ is the set of 
real-valued predictable processes $\phi$ such that $$\| \phi\|^p_{\H^{p,T}} := E \left[(\int_0 ^T \phi_t ^2 dt)^{\frac{p}{2}}\right] < \infty.$$  

For $\beta >0$ and  $\phi \in \H^{2,T} $, we introduce the norm $\| \phi \|_{\beta,T}^2 := E[\int_0^T e^{\beta s} \phi_s^2 ds]. $

\item $L^p_\nu$ is the set of Borelian functions $\ell: \R^* \rightarrow \R$ such that  $\int_{\R^*}  |\ell(u) |^p \nu(du) < + \infty.$

The set  $L^2_\nu$ is a 
 Hilbert space equipped with the scalar product 
$$\langle \delta    , \, \ell \rangle_\nu := \int_{\R^*} \delta(u) \ell(u) \nu(du)  \quad \text{ for all  }  \delta  , \, \ell \in L^2_\nu \times L^2_\nu,$$ 
and the norm $\|\ell\|_\nu^2 :=\int_{\R^*}  |\ell(u) |^2 \nu(du) < + \infty.$


\item $\H_{\nu}^{p,T}$ is  the set of processes $l$ which are {\em predictable}, that is, measurable $$l : ([0,T]  \times \Omega \times \R^*,\; \PC \otimes {\cal B}(\R^*))  \rightarrow (\R\;,  \BC(\R)); \quad 
(\omega,t,u) \mapsto l_t(\omega, u)
$$
such that $$\| l \|^p_{\H_{\nu}^{p,T}} :=E\left[( \int_0 ^T \|l_t\|_{\nu}^2 \,dt ) ^{\frac{p}{2}} \right]< \infty.$$

For $\beta >0$ and $l \in \H_\nu^{2,T}$, we set 
 $\| l \|_{\nu,\beta,T}^2 := E[\int_0^T e^{\beta s} \|l_s\|_\nu^2 \, ds] $. 

\item   ${\cal S}^{p,T}$ is the set of real-valued RCLL adapted 
 processes $\phi$ such that  
$$\| \phi\|^p_{{\cal S}^p} := E(\sup_{0\leq t \leq T} |\phi_t |^p) <  \infty.$$ 



\end {itemize}

When $T$ is fixed and there is no ambiguity, we denote
$\H^{p}$ instead of $\H^{p,T}$,   $\H_{\nu}^{p}$ instead of $\H_{\nu}^{p,T}$, 
${\cal S}^p$ instead of  ${\cal S}^{p,T}$.

\begin{itemize}

\item  $\T_{0}$ denotes the set of 
stopping times $\tau$ such that $\tau \in [0,T]$ a.s.

\item For $S$ in $\T_{0}$,    $\T_{S}$  is the set of 
stopping times 
$\tau$ such that $S \leq \tau \leq T$ a.s.

\end {itemize}

\begin{definition}[Driver, Lipschitz driver]\label{defd}
A function $f$ is said to be a {\em driver} if 
\begin{itemize}
\item  
$f: [0,T]  \times \Omega \times \R^2 \times L^2_\nu \rightarrow \R $\\
$(\omega, t,x, \pi, \ell(\cdot)) \mapsto  f(\omega, t,x, \pi, \ell(\cdot)) $
  is $ {\cal P} \otimes {\cal B}(\R^2)  \otimes {\cal B}(L^2_\nu) 
- $ measurable,  
\item $f(.,0,0,0) \in \H^2$.
\end{itemize} 
A driver $f$ is called a {\em Lipschitz driver} if moreover there exists a constant $ C \geq 0$ such that $dP \otimes dt$-a.s.\,, 
for each $(x_1, \pi_1, \ell_1)$, $(x_2, \pi_2, \ell_2)$, 
$$|f(\omega, t, x_1, \pi_1, \ell_1) - f(\omega, t, x_2, \pi_2, \ell_2)| \leq 
C (|x_1 - x_2| + |\pi_1 - \pi_2| +   \|\ell_1 - \ell_2 \|_\nu).$$
\end{definition}


%

%
\paragraph{Existence and uniqueness result for BSDEs with jumps.} (Tang and Li ,1994  \cite{Tang})
 Let $T >0$. For each Lipschitz driver $f$, and each terminal condition $\xi$ $\in$ $L^2({\cal F}_T)$, there exists a unique solution $(X, \pi, l)$ $\in$ ${\cal S}^{2,T} \times \H^{2,T} \times \H_{\nu}^{2,T}$ satisfying
 \begin{equation}\label{er}
-dX_t =f(t,X_{t^-},\pi_t, l_t(\cdot))dt - \pi_t dW_t -  \int_{\R^*} l_t(u) \tilde N(dt,du)  ; \qquad  X_T = \xi. 
\end{equation} 
 This  solution is denoted by $(X( \xi,T), \pi ( \xi,T), l( \xi,T))$.

This result can be extended if the terminal time  $T$ is replaced by a 
stopping time $S \in \T_0$.  Let $(X(\xi,S), \pi(\xi,S), l(\xi,S))$ (denoted here by $(X,\pi,l)$) be the solution of the BSDE associated with driver $f$, terminal time $S$ and terminal condition $\xi \in L^2({\cal F}_S)$. The solution can be extended on the whole interval $[0,T]$ by setting $X_t = \xi, \pi_t =0,l_t=0$ for $t \geq S$. So, $\left((X_t, \pi_t, l_t); t \leq T\right)$ is the unique solution of the BSDE with  driver $f(t,x, \pi,l) {\bf 1}_ { \{ t \leq S  \}  }$ and terminal conditions ($T$, $\xi$).

We  refer to  \cite{BBP,R} and to  \cite{QuenSul} where some  results are used in this paper.
%
%


\paragraph{Dynamic risk measures induced by BSDEs with jumps.}

Let $T' >0$ be a time horizon. Let $f$ be a Lipschitz driver such that $f(\cdot,0,0,0) \in \H^{2,T'}$.
We define   the following functional:
 for each $T \in [0,T']$ and $\xi \in L^2({\cal F}_T)$, set 
  \begin{equation}\label{definition}
  \rho_t^f (\xi, T) = \rho_t (\xi, T) :=  -X_t(\xi,T), \,\, \,\,\,0\leq t \leq T.
  \end{equation}
 where $X_t(\xi,T)$ denotes the solution  of the BSDE \eqref{er}  with driver $f$, terminal condition $\xi$ and terminal time  $T$.
If $T$ represents a given maturity and $\xi$ a financial position at time $T$,
then $\rho_t (\xi, T)$ will be interpreted as the risk of $\xi$ at time 
$t$. The functional $\rho :  (\xi , T) \mapsto \rho_\cdot(\xi, T) $ defines then
 a dynamic risk measure induced by the BSDE with driver $f$.
Properties of such dynamic risk measures are given in \cite{QuenSul}.

\paragraph{Optimal stopping problem.}
The aim of this paper is to study optimal stopping for  dynamic risk measures. 
 Let $T >0$ be the terminal time.
 Let $\{\xi_t, 0 \leq t \leq T \}$ be a RCLL 
 adapted process on $[0, T]$, belonging to ${\cal S}^2$, representing a dynamic financial position. 
  

Consider the following optimal stopping problem:
For each stopping time $S$ $\in$ $\T_0$, let $v(S)$ be the ${\cal F}_S$-measurable random variable (unique for the equality in the almost sure sense) defined by 
 \begin{eqnarray}\label{un}
v(S):= {\rm ess} \inf_{\tau \in \T_S}\rho_S(\xi_{\tau},\tau).
\end{eqnarray}
Since by definition $\rho_S(\xi_{\tau},\tau)= -X_S(\xi_{\tau}, \tau)$, we  have that for each 
stopping time $S \in \T_0$, 
\begin{equation} \label{vvv}
v(S) =  {\rm ess} \inf_{\tau \in \T_S}-X_S(\xi_{\tau}, \tau)=-  {\rm ess} \sup_{\tau \in \T_S}X_S(\xi_{\tau}, \tau).
\end{equation}
The aim is to characterize  for each $S \in \T_S$ the minimal risk-measure $v(S)$ and  to provide an existence 
result of 
an $S$-optimal stopping time $\tau^* \in \T_S$, that is
such that   $v(S) = \rho_S(\xi_{\tau^*}, \tau^*)$ a.s.  
This problem is related to reflected BSDEs. 
We give below existence and uniqueness results for these equations. 

\section{RBSDEs with jumps and RCLL obstacle process}\label{rbsde}
Reflected BSDEs  (RBSDEs) have been introduced by  N. El Karoui et al. (1997 ) (see \cite{ElKaroui97}). The solution of such equations 
are constrained to be greater than a given process called the obstacle.
In this section, we provide existence and uniqueness results for RBSDEs with jumps, in the case when the obstacle is RCLL,  which  complete some results in \cite{HO1, HO2, Essaky}.

Let $T>0$ be a fixed terminal time and  $f$ be  a Lipschitz driver. 
Let $\xi_.$ be a process called obstacle in ${\cal S}^2 $. 

\begin{definition}
A process $(Y,Z,k(.),A)$ is said to be a solution of the reflected BSDE associated with driver $f$ and obstacle $\xi_.$ if 
\begin{align}\label{RBSDE}
&(Y,Z,k(.),A)\in {\cal S}^2 \times \H^2 \times \H^2_{\nu}\times {\cal S}^2 \nonumber \\
& -dY_t = f(t,Y_t,  Z_t, k_t(\cdot) )dt +dA_t - Z_t  dW_t - \int_{\R^*} k_t(u) \tilde{N}(dt,du); \quad  Y_T = \xi_T, \\
&  Y_t \geq \xi_t,  \; 0 \leq t \leq T \text{ a.s.}, \nonumber\\
& A \text{ is a nondecreasing RCLL predictable process 
with } A_0= 0 \text{ and such that } \nonumber\\
& \int_0^T (Y_t - \xi_t) dA^c_t = 0 \text{ a.s. and } \;  \Delta A_t^d= - \Delta Y_t {\bf 1}_{Y_{t^-} = \xi_{t^-}} \text{ a.s.} \nonumber
\end{align}
\end{definition}
Here  $A^c$ denotes the continuous part of $A$ and $A^d$  its discontinuous part.

We  introduce the following definition.

\begin{Definition} \label{defr} A progressive process $(\phi_t)$ is said to be {\em left-upper semicontinuous along stopping times} if for all $\tau \in \T_0$ and for each non decreasing sequence of stopping times $ (\tau_n)$ such that $\tau^n \uparrow \tau$ a.s.\,,
\begin{equation}\label{usc}
\phi_{\tau} \geq \limsup_{n\to \infty} \phi_{\tau_n} \quad \mbox{a.s.}
\end{equation} 
\end{Definition}

\begin{remark}
Note that in this definition, no condition is required at a totally unaccessible
stopping time. In our framework,  since the filtration is generated by $W$ and $N$, this means that  no condition is required at the jump times of $N$.
\end{remark}

\subsection{The case when the driver $f$ does not depend on $y,z,k$.}

\begin{Proposition}\label{f}
Suppose that $f$ does not depend on $y,z,k$, that is $f(\omega,t, y,z, k(\cdot)) = f(\omega,t)$, where $f$ 
is in $\H^2 := \H^{2,T}$. 
Then, RBSDE~\eqref{RBSDE} admits a unique solution $(Y,Z,k(.),A)\in {\cal S}^2 \times \H^2 \times \H^2_{\nu}\times {\cal S}^2$ and  for each $S \in \T_0$,  
\begin{eqnarray}\label{deux-2}
 Y_S= {\rm ess} \sup_{\tau \in \T_S} E[ \xi_{\tau} + \int_\tau^T f(t)dt \mid \Fc_S] \quad \mbox{a.s.}
\end{eqnarray} 
Moreover if $(\xi_t)$ is  left-upper semicontinuous along stopping times, then $A_t$ is continuous.
\end{Proposition}
\dproof
For each $S \in \T_0$,  we introduce the following random variable 
\begin{eqnarray}\label{ddeux-2}
\overline Y  (S):= {\rm ess} \sup_{\tau \in \T_S} E[ \xi_{\tau} + \int_\tau^T f(t)dt \mid \Fc_S]. 
\end{eqnarray} 
By classical results of optimal control theory, there exists a RCLL adapted process denoted by $(\overline Y_t)$ such that 
for each $S \in \T_0$, 
$\overline{Y}(S) = \overline Y_S$ a.s. 
The process $(\overline Y_t + \int_0^t f(s) ds)$ is a supermartingale. By the Doob-Meyer decomposition, it can be uniquely written as 
$$
d\overline Y_t = -f(t)dt -dA_t + dM_t,
$$
where $M$ is a square-integrable martingale and $A$ is a nondecreasing RCLL predictable process with 
$E (A_T^2) < \infty$  and  $A_0= 0$. 
Furthermore, by the theorem of representation \cite{Tang}, there exist unique processes $Z$ in $\H^2$ and $k$ in $\H^2_\nu$ such that 
$$dM_t = Z_t  dW_t + \int_{\R^*} k_t(u) \tilde{N}(dt,du).$$
The process $A$ can be uniquely decomposed as  $dA_t = dA^c_t + dA_t^d$. 
By Proposition B.11 in \cite{KQ} (or \cite{EK}), we have
$\int_0^T (\overline Y_t - \xi_t) dA^c_t = 0$ a.s. and $\Delta A_t^d= - \Delta Y_t 
{\bf 1}_{\overline Y_{t^-} = \xi_{t^-}}$ a.s. Hence, $(\overline Y, Z, k(), A)$ is a solution of the   RBSDE associated with driver $f(t)$ and obstacle $(\xi_t)$.

In the particular case when $(\xi_t)$ is left-upper semicontinuous over stopping times, by Proposition 2.11 in \cite{KQ} (see also \cite{EK}), the supermartingale $v_t:=\overline Y_t + \int_0^t f(s) ds$ is then left-continuous over stopping times in expectation, that is, for all $\tau \in \T_0$ and for each non decreasing sequence of stopping times $ (\tau_n)$ such that $\tau^n \uparrow \tau$ a.s.\,, $\lim_{n\to \infty} E[v_{\tau_n}]= E[v_{\tau}]$. Consequently, by Lem. B.8 in \cite{KQ} (or Th. 10, Chap. VII in \cite{DM2}),
the nondecreasing process $A$ is continuous.

We will now show that conversely, if $(Y, Z, k(\cdot), A)$ is a solution of the RBSDE associated with driver $f(t)$ and obstacle $(\xi_t)$, then, for each $S$ $\in$ $\T_0$, $Y_S = \overline Y(S)$ a.s.\, 

To simplify, suppose that $f= 0$. The following proof can be easily generalized to the case where $f\neq 0$. Suppose that $(Y,Z,k(.),A)$ is a solution of the reflected BSDE associated with driver $f=0$ and obstacle $\xi_t$. For each $t$, let 
$$M_t:= Y_0 +\int_0^t Z_s  dW_s +\int_0^t  \int_{\R^*} k_s(u) \tilde{N}(ds,du).$$
Note that $M$ is a square integrable martingale. We have
\begin{eqnarray}\label{RBSDE1}
dY_t = dM_t - {\bf 1}_{\{ Y_t = \xi_t  \} } dA^c_t - dA^d_t;  && 
 Y_T = \xi_T, 
\end{eqnarray}
with $\Delta A_t^d= - \Delta Y_t 
{\bf 1}_{Y_{t^-} = \xi_{t^-}}$ a.s. Since $Y \geq \xi$, it clearly follows that for each stopping time $S$ $\in$ $\T_0$ and for each $\tau$ $\in$ $\T_S,$
\begin{eqnarray*}
Y_S= E[ Y_{\tau}  \mid \Fc_S]-  E[ A_{\tau}- A_S  \mid \Fc_S]  \geq E[ \xi_{\tau}  \mid \Fc_S] \quad \mbox{a.s.}
\end{eqnarray*} 
Hence, by taking the supremum over $\tau$ $\in$ $\T_S,$ we have 
 \begin{eqnarray}\label{firsti1}
Y_S \geq {\rm ess} \sup_{\tau \in \T_S}E[ \xi_{\tau}  \mid \Fc_S]= \overline Y(S) \quad \mbox{a.s.}
\end{eqnarray}
It remains to show the converse inequality.

Let us first consider the simpler case where $(\xi_t)$ is  left-upper semicontinuous over stopping times  and $A$ is continuous, that is $A=A^c$. For each $S$ $\in$ $\T_0$, consider 
$$\tau^*_S:= \inf \{ t \geq S,\,\, Y_t = \xi_t\}.$$
Note that $\tau^*_S$ $\in$ $\T_S$. Since $Y$ and $\xi$ are right-continuous processes, we have $Y_ { \tau^*_S } = \xi_{  \tau^*_S}$ a.s. 
By definition of $\tau^*_S$, for almost every $\omega$, for each $t \in$ $[S(\omega), \tau^*_S(\omega)[$, we have $Y_t(\omega)> \xi_t(\omega)$. Hence, since $Y$ is solution of the RBSDE, for almost every $\omega$, the nondecreasing function $t \mapsto A_t(\omega)$ is constant on $[S(\omega), \tau^*_S(\omega)[$. The continuity of $A$ implies that $t \mapsto A_t(\omega)$ is constant on $[S(\omega), \tau^*_S(\omega)]$. This clearly leads to the following equality:
$$Y_S= E[ \xi_{\tau^*_S}  \mid \Fc_S] \quad \mbox{a.s.}$$
This with inequality (\ref{firsti1}) gives the desired equality $Y_S = \overline Y(S)$ a.s.

We now consider the case where $(\xi_t)$ is only supposed to be a RCLL process. For each $S$ $\in$ $\T_0$ and for each $\varepsilon >0$, let 
$$\tau^{\varepsilon}_S:= \inf \{ t \geq S,\,\,Y_t \leq \xi_t  + \varepsilon\}.$$
Note that $\tau^{\varepsilon}_S$ $\in$ $\T_S$. Fix $\varepsilon >0$. For a.e. $\omega$, if $t \in [S(\omega), 
\tau^{\varepsilon}_S(\omega)[$, then $ Y_t(\omega)> \xi_t (\omega) + \varepsilon$ and hence $Y_t(\omega)> \xi_t (\omega)$. It follows that for a.e. $\omega$, the function $t \mapsto A^c_t (\omega)$ is constant on $[S(\omega), 
\tau^{\varepsilon}_S(\omega)]$ and $t \mapsto A^d_t (\omega)$ is constant on $[S(\omega), \tau^{\varepsilon}_S(\omega)[$.
Also, $ Y_{ (\tau^{\varepsilon}_S)^-  } \geq \xi_{ (\tau^{\varepsilon}_S)^-  }+ \varepsilon\,$ a.s.\,
Since $\varepsilon >0$, it follows that $ Y_{ (\tau^{\varepsilon}_S)^-  } > \xi_{ (\tau^{\varepsilon}_S)^-  }$  a.s.\,\,, which implies that $\Delta A^d _ {\tau^{\varepsilon}_S  } =0$ a.s.\, The process $(Y_t)$ is thus a martingale on $[S, 
\tau^{\varepsilon}_S]$. Furthermore, by the right-continuity of $(\xi_t)$ and $(Y_t)$, we clearly have 
$$Y_{\tau^{\varepsilon}_S} \leq \xi_{\tau^{\varepsilon}_S} + \varepsilon \quad \mbox{a.s.}$$
It follows that 
\begin{eqnarray}\label{lambda1}
 Y_S =   E[ Y_{\tau^{\varepsilon}_S}  \mid \Fc_S] \leq E[ \xi_{\tau^{\varepsilon}_S}  \mid \Fc_S]  + \varepsilon \leq  {\rm ess} \sup_{\tau \in \T_S}E[ \xi_{\tau}  \mid \Fc_S] + \varepsilon \quad \mbox{a.s.}
\end{eqnarray}
Hence, 
$Y_S \leq \overline Y(S)  + \varepsilon$ a.s.\,
for each $\varepsilon >0$, which implies that $Y_S \leq  \overline Y(S)$ a.s.\, This, with inequality \eqref{firsti1}, ensures the desired equality $Y_S = \overline Y(S)$ a.s. 
\fproof


\subsection{The case of a general Lipschitz  driver}
\begin{theorem}\label{exiuni}
Suppose that $f$ is a  Lipschitz driver with Lipschitz constant $C$. Then, RBSDE (\ref{RBSDE}) admits a unique solution $(Y,Z,k(.),A)\in {\cal S}^2 \times \H^2 \times \H^2_{\nu}\times {\cal S}^2$. Moreover if $(\xi_t)$ is  left-upper semicontinuous over stopping times, then $A_t$ is continuous.
\end{theorem}
\dproof
We denote by $\H_\beta^2$ the space $\H^2 \times \H^2 \times \H^2_\nu$ equipped with the norm 
$\| Y, Z, k(\cdot) \|_\beta^2 := \| Y \|_{\beta}^2 +  \| Z \|_{\beta}^2  +  \| k \|_{\nu,\beta}^2 $.

We define a mapping $\Phi$ from $\H_\beta^2$ into
itself as follows. Given $(U, V,l) \in \H_\beta^2$, let $(Y, Z,k) = \Phi (U, V, l)$
be the 
the solution of the RBSDE  associated with driver $f(s) = 
f(s, U_s , V_s, l_s)$.  
Let $A$ be the associated nondecreasing process. 
The mapping $\Phi$ is well defined by Proposition \ref{f}.
By using some a priori estimates (see Proposition \ref{est}), $\Phi$ can be shown to be a contraction from $\H_\beta^2$ into
itself. It thus admits an unique fixed point, which corresponds to the solution of RBSDE (\ref{RBSDE}). 
For details, see the Appendix.
\fproof

\section{Relations between optimal stopping problems and   RBSDEs }\label{sec-charact}

In the following, we make the following assumption on the driver $f$  which ensures  the monotonicity property of the associated risk measure  $\rho$   (see \cite{QuenSul}).

Let $T >0$.
\begin{Assumption}\label{Royer}
A driver $f$ is said to satisfy Assumption~\ref{Royer} if the following holds:\\
$dP \otimes dt$-a.s\, for each $(x,\pi,l_1,l_2)$ $\in$ $[0,T] \times \Omega \times \RB^2 \times (L^2_{\nu})^2$,
$$f( t,x,\pi,l_1)- f(t,x,\pi,l_2) \geq \langle \theta_t^{x,\pi,l_1,l_2}  \,,\,l_1 - l_2 \rangle_\nu,$$ 
with
\begin{equation*}
\theta:  [0,T]  \times \Omega\times \RB^2 \times  (L^2_{\nu})^2  \mapsto  L^2_{\nu}\,; \, (\omega, t, x,\pi,l_1,l_2) \mapsto 
\theta_t^{x,\pi,l_1,l_2}(\omega,.)
\end{equation*}
 ${\cal P } ·\otimes {\cal B}(\R^2) \otimes  {\cal B}( (L^2_{\nu})^2 )$-measurable, bounded, and satisfying $ dP\otimes dt \otimes d\nu(u)$-a.s.\,, for each $(x,\pi,l_1,l_2)$ $\in$ $\RB^2 \times (L^2_{\nu})^2$,
 \begin{equation}\label{condi}
\theta_t^{x,\pi,l_1,l_2} (u)\geq -1\,\,\,  \;\; \text{ and }
\,\,  \;\;|\theta_t^{x,\pi,l_1,l_2}(u)|  \leq \psi(u),
\end{equation}
where $\psi$ $\in$ $L^2_{\nu}$.
\end{Assumption}

\subsection{Characterization of  the value function  as the solution of  an RBSDE}

We relate the  optimal stopping problem \eqref{vvv} to reflected BSDEs. We first   show that  the value function $v$ coincides with $-Y$, where $Y$ is the solution of the reflected BSDE associated with driver $f$ and obstacle $\xi$. 

\begin{theorem}[Characterization]\label{caracterisation}
Let $T >0$ be the terminal time.
 Let $(\xi_t, 0 
\leq t \leq T )$ be a RCLL 
 adapted process on $[0, T]$, belonging to ${\cal S}^2$.
Suppose that $(Y,Z,k(\cdot),A)$ is the solution of the reflected BSDE (\ref{RBSDE}).
Then, for each stopping time $S$ $\in$ $\T_0$, we have
 \begin{equation}\label{prixam}
Y_S =  {\rm ess} \sup_{\tau \in \T_S}X_S(\xi_{\tau}, \tau) \quad \mbox{a.s.}
\end{equation}
where on the interval 
$[S, \tau]$, the process $(X_{s}(\tau, \xi_{\tau}),
\pi_{s}(\tau, \xi_{\tau}), l_s(\tau, \xi_{\tau}))$ satisfies the BSDE
\begin{eqnarray*} 
 -dX_{s}= f(s,X_{s}, \pi_{s}, l_s)ds - \pi_{s}dW_{s} - \int_{\R^*} l_s(u) \tilde{N}(ds,du); 
&& X_{\tau}= \xi_{\tau}.
 \end{eqnarray*}

\end{theorem}

\dproof  

We first show that $Y_{S} \geq X_S(\tau, \xi_{\tau})$, for 
each
$\tau$ $\in$ $\T_0$.  
Fix $\tau$ $\in$ $\T_0$.  
In the interval $[S,\tau]$, the process $(Y,Z,k(\cdot),A)$
 satisfies~:
 \begin{eqnarray*} 
 -dY_{s}= f(s,Y_{s}, Z_{s}, k_s)ds + dA_s - Z_{s}dW_{s} - \int_{\R^*} k_s(u) \tilde{N}(ds,du); 
&& Y_{\tau}= Y_{\tau}.
 \end{eqnarray*}
In other words,  the process $(Y_{s}, Z_{s}, k_s;  S \leq s \leq \tau)$ is the 
solution 
 of the  BSDE associated with terminal time $\tau$, terminal
condition $Y_{\tau}$ and (generalized) driver 
$$f(s,y,z ,k)ds+ dA_{s}.$$
Since $f(s,y,z ,k)ds + dA_{s} \geq f(s,y,z,k)ds $ and since $Y_{\tau} \geq
\xi_{\tau}$ a.s.\,, the comparison theorem for BSDEs (see Theorem 4.2 in  \cite{QuenSul})  gives that
$$Y_{S} \geq X_S( \xi_{\tau}, \tau) \quad \mbox{a.s.}$$
By taking the supremum over $\tau \in \T_S$, we derive that 
 \begin{equation}\label{prixame}
Y_S \geq  {\rm ess} \sup_{\tau \in \T_S}X_S(\xi_{\tau}, \tau) \quad \mbox{a.s.}
\end{equation}

It remains to show the converse inequality.\\ 
For each $S$ $\in$ $\T_0$ and for each $\varepsilon >0$, let 
$\tau^{\varepsilon}_S$ be the stopping time defined by
\begin{equation}\label{thetalambda}
\tau^{\varepsilon}_S := \inf \{ t \geq S,\,\, Y_t \leq \xi_t + \varepsilon\}.
\end{equation}

We first show two useful lemmas.
\begin{lemma}\label{lala}  
\begin{itemize}
\item We have
\begin{equation*}
 Y_{\theta^{\varepsilon}_S} \leq \xi_{\theta^{\varepsilon}_S} + \varepsilon \quad \mbox{a.s.}
\end{equation*}
\item
The process 
$(Y_t, S\leq t \leq 
\theta^{\varepsilon}_S )$ is the solution of the BSDE associated with terminal time $\theta^{\varepsilon}_S $, terminal
condition $Y_{\theta^{\varepsilon}_S }$ and driver $f$, that is
\begin{equation*}
Y_t = X_t (Y_{\theta^{\varepsilon}_S }, \theta^{\varepsilon}_S) \quad S \leq t \leq 
\theta^{\varepsilon}_S \ 
 \mbox{ a.s.}
\end{equation*}
\end{itemize}
\end{lemma}

\dproof
The first point follows from the definition of $\theta^{\varepsilon}_S$ and the right-continuity of $(\xi_t)$ and $(Y_t)$. Let us show the second point.
Note that $\theta^{\varepsilon}_S$ $\in$ $\T_S$. Fix $\varepsilon >0$. 
For a.e. $\omega$, if $t \in [S(\omega), 
\theta^{\varepsilon}_S(\omega)[$, then $Y_t(\omega)> \xi_t (\omega)+ \varepsilon$ and hence $Y_t(\omega)> \xi_t (\omega)$. It follows that for a.e. $\omega$, the function $t \mapsto A^c_t (\omega)$ is constant on $[S(\omega), 
\theta^{\varepsilon}_S(\omega)]$ and $t \mapsto A^d_t (\omega)$ is constant on $[S(\omega), \theta^{\varepsilon}_S(\omega)[$.
Also, $ Y_{ (\theta^{\varepsilon}_S)^-  } \geq \xi_{ (\theta^{\varepsilon}_S)^-  }+ \varepsilon \,$ a.s.\,\,
Since $\varepsilon >0$, it follows that $ Y_{ (\theta^{\varepsilon}_S)^-  } > \xi_{ (\theta^{\varepsilon}_S)^-  }\quad \mbox{a.s.}$, which implies that $\Delta A^d _ {\theta^{\varepsilon}_S  } =0$ a.s.\, 
\fproof

\begin{lemma}\label{eps}
Set $\beta:= 3 C^2 +2C$, where $C$ is the Lipschitz constant of $f$.\\
For each $\varepsilon$ $>0$ and each $S$ $\in$ $\T_0$, we have
\begin{equation}\label{fifi}
Y_S \,\,\leq \,\,  X_S ( \xi_{\theta^{\varepsilon}_S}, \theta^{\varepsilon}_S)  + e^{\frac{\beta T} {2}}\varepsilon \quad \mbox{a.s.}\,
\end{equation}

 \end{lemma}
 
 \dproof
By Lemma \ref{lala} and by the comparison theorem for BSDEs, we derive that for each $\varepsilon >0$,
\begin{equation}\label{fi}
  Y_S = X_S (Y_{\theta^{\varepsilon}_S }, \theta^{\varepsilon}_S)  \leq  X_S(\xi_{\theta^{\varepsilon}_S}+ \varepsilon, \theta^{\varepsilon}_S)  \quad \mbox{a.s.}
\end{equation}
Now, by the a priori estimates on BSDEs (see Proposition A.4 \cite{QuenSul}), we have
\begin{equation*}
 \vert X_S(\xi_{\theta^{\varepsilon}_S}+ \varepsilon, \theta^{\varepsilon}_S) 
- X_S ( \xi_{\theta^{\varepsilon}_S}, \theta^{\varepsilon}_S) \vert^2
\leq e^{\beta (T-S)} \varepsilon^2\,\quad \mbox{a.s.}
\end{equation*}
This with inequality (\ref{fi}) leads to inequality (\ref{fifi}), which ends the proof of Lemma \ref{eps}.
\fproof
 
 {\bf End of the proof of Theorem \ref{caracterisation}}\\
By Lemma \ref{eps}, we have for each $\varepsilon >0$,
\begin{equation}\label{etoile}
Y_S \,\,\leq \,\,  X_S ( \xi_{\theta^{\varepsilon}_S}, \theta^{\varepsilon}_S)  + e^{\frac{\beta T} {2}}\varepsilon
\,\, \leq  \,\, {\rm ess} \sup_{\tau \in \T_S}X_S(\xi_{\tau}, \tau) + e^{\frac{\beta T} {2}  } \varepsilon  \quad \mbox{a.s.}
\end{equation}
 It follows that 
 \begin{eqnarray*}
Y_S \,\, \leq  \,\,{\rm ess} \sup_{\tau \in \T_S}X_S(\xi_{\tau}, \tau) \quad \mbox{a.s.}\,,
\end{eqnarray*}
and, since we have already shown the converse inequality, this inequality is an equality.
\fproof

\begin{remark}
By inequality \eqref {fifi}, the stopping time $\theta^{\varepsilon}_S$ is an $\varepsilon '$-optimal stopping time for the optimal stopping time problem \eqref{prixam} with $\varepsilon ' =  e^{\frac{\beta T} {2}}\varepsilon$.

Note also that the above result does not require any concavity assumption on the driver, contrary to 
 \cite{ Bayraktar}  and \cite{ Bayraktar-2}.

\end{remark}

  

\subsection{Optimal stopping times }

We now provide an optimality criterium  for the optimal stopping time problem \eqref{prixam}.
\begin{Proposition}[Optimality criterium.] \label{optcri}
Let $S \in \T_0$ and let $\hat{\tau} \in  \T_S.$ Suppose that in Assumption \ref{Royer}
 for each $l$ $\in$ ${\cal L}^2_\nu$, we have
 \begin{equation}\label{assu}
\theta_t^{X^{\hat \tau},\pi^{\hat \tau},l\,, l^{\hat \tau}} > - 1, \quad dt \otimes dP-{\rm a.s.}\,
\end{equation}
where $(X^{\hat \tau},\pi^{\hat \tau}, l^{\hat \tau}):= 
(X(\xi_{\hat \tau}, \hat{\tau}), \pi(\xi_{\hat \tau}, \hat{\tau}), l(\xi_{\hat \tau}, \hat{\tau}))$ is  the solution of the BSDE associated with $\hat\tau$, $\xi_{\hat \tau}$.\\
The stopping time $\hat{\tau}$ is $S$-optimal, i.e.
\begin{equation}\label{uni}
Y_S =  {\rm ess} \sup_{\tau \in \T_S}X_S(\xi_{\tau}, \tau)=X_S(\xi_{\hat \tau}, \hat{\tau}) \; \; \text{ a.s. } 
\end{equation}
if and only if 
\begin{equation}\label{deuxi}Y_s = X_s(\xi_{\hat \tau}, \hat{\tau}), \quad S \leq s \leq 
\hat{\tau}  \text{ a .s. } 
\end{equation}
In other words,  $\hat{\tau}$ is $S$-optimal if and only if $(Y_s,  \,S \leq s \leq  {\hat \tau})$ is the solution of the non reflected BSDE associated with terminal time ${\hat \tau}$  and terminal condition $\xi_{\hat \tau}$.
\end{Proposition}

\dproof
It is clear that (\ref{deuxi}) $\Rightarrow $ (\ref{uni}). Note that this implication does not require condition (\ref{assu}).
It remains to prove that (\ref{uni}) $\Rightarrow $ (\ref{deuxi}).

Suppose that $\hat{\tau}$ is an $S$-optimal stopping time. 

The process $(Y_{s}, Z_{s}, k_s;  S \leq s \leq \hat{\tau})$ is the 
solution 
 of the  BSDE associated with terminal time $\hat{\tau}$, terminal
condition $Y_{\hat{\tau}}$ and (generalized) driver 
$$f(s,y,z,k)ds + dA_{s}.$$
We have $f(s,y,z ,k)ds + dA_{s} \geq f(s,y,z,k)ds $, 
 $Y_{\hat{\tau}} \geq
\xi_{\hat{\tau}}$ a.s. as well as equality \eqref{uni}. 
Using Assumption (\ref{assu}) and  applying the strict comparison theorem for BSDEs (see \cite{QuenSul}, Th 4.4), we get the desired result.

 \fproof
 \begin{remark}
In the particular case when the driver $f$ does not depend on $(y,z)$, this gives  the well-known optimality criterium of the Optimal Stopping Theory: a stopping time $\hat{\tau}$ is $S$-optimal if and only if $(Y_s + \int_0^s f(r)dr),  \,S \leq s \leq  {\hat \tau})$ is a martingale with $Y_{\hat \tau}=\xi_{\hat \tau}$ a.s.
\end{remark}

We now show that,  under a left regularity condition on the obstacle,  
$\tau^{\varepsilon}_S$ tends to an $S$- optimal stopping time for Problem~\eqref{prixam} as $\varepsilon$ tends to $0$\,,
and we provide some additional properties.

\begin{theorem}\label{existence}
Suppose $(\xi_t)$ is left-upper semi-continuous along stopping times.
Let $S \in \T_0$.

(i) The stopping time $\tilde \tau_S$ defined by 
$${\tilde \tau}_S := \lim_{\varepsilon \downarrow 0}\uparrow \tau^{\varepsilon}_S$$
is an $S$-optimal stopping time.

(ii) the stopping time $\tau_S^*$ defined by 
$$\tau^*_S := \inf \{u \geq S ; \, Y_u = \xi_u \}$$ 
is an $S$-optimal stopping time and we have
$$Y_s = X_s(\xi_{\tau^*_S}, \tau^*_S), \quad S \leq s \leq 
\tau^*_S  \text{ a .s. } $$
We also have $\tau_S^* \geq \tilde \tau_S$ a.s.

(iii) Suppose moreover that in Assumption \ref{Royer},  for all $x, \pi,l_1,l_2$, we have
\begin{equation}\label{hypsup}
\theta_t^{x,\pi,l_1,l_2} > -1 \quad dt \otimes dP-{\rm a.s.}
\end{equation}
Then,  $\tau_S^* = \tilde \tau_S$ a.s. and $\tau_S^*$ is the minimal $S$-optimal stopping time.
\end{theorem}

\dproof

(i) By letting $\varepsilon$  tend to $0$ in inequality \eqref{etoile}, we get
\begin{equation}\label{lamb}
 Y_S   \leq   \limsup_{\varepsilon \downarrow 0} X_S(\xi_{\tau^{\varepsilon}_S}, \tau^{\varepsilon}_S).
\end{equation}
For each $\omega$ such that the map $\varepsilon \mapsto \tau^{\varepsilon}_S(\omega)$ from $\R^*_+ \rightarrow [0,T]$ is constant for $\varepsilon$ sufficiently small, we have
 $$\lim_{\varepsilon \downarrow 0} \xi_{\tau^{\varepsilon}_S}(\omega) = \xi_{{\tilde \tau}_{S}}(\omega).$$
Moreover, since the process $(\xi_t)$ is left-limited,  for almost every $\omega$ such that for each  $\varepsilon >0$, $\tau^{\varepsilon}_S(\omega)  < 
\hat{\tau}_{S}(\omega)$, we have
$$\lim_{\varepsilon \downarrow 0} \xi_{\tau^{\varepsilon}_S}(\omega) = \xi_{{\tilde \tau}_{S}^-}(\omega).$$
Hence, for almost every $\omega$, $\lim_{\varepsilon \downarrow 0} \xi_{\tau^{\varepsilon}_S}(\omega)$ does exist.

The continuity property of BSDEs
with respect to terminal conditions (see Prop. A6 in \cite{QuenSul}), 
implies
\begin{equation}\label{lamb2}
\lim_{\varepsilon \downarrow 0} X_S(\xi_{\tau^{\varepsilon}_S}, \tau^{\varepsilon}_S) 
=  X_S(\lim_{\varepsilon \downarrow 0} \xi_{\tau^{\varepsilon}_S},\tilde {\tau}_S)
 \quad \mbox{a.s.}
\end{equation}
Now, by the left-upper semicontinuity property of 
the obstacle along stopping times, we have
$$\lim_{\varepsilon \downarrow 0}\xi_{\tau^{\varepsilon}_S} \leq \xi_{\tilde{\tau}_S} 
\; \text{ a.s. }$$
By the comparison theorem, it follows that
$$ X_S(\lim_{\varepsilon \downarrow 0} \xi_{\tau^{\varepsilon}_S},\tilde {\tau}_S)
 \leq X_S(\xi_{{\tilde \tau}_{S}}, \tilde \tau_S).$$
Hence, by \eqref{lamb} and \eqref{lamb2}, we get $Y_S   \leq X_S(\xi_{{\tilde \tau}_{S}}, \tilde \tau_S)$ a.s.\, By using the characterization of $Y_S$ as the value function of the optimal stopping time problem~\eqref{prixam}, we get
\begin{equation}  Y_S   = X_S(\xi_{{\tilde \tau}_{S}}, \tilde \tau_S) \quad \mbox{a.s.}
\end{equation}
Thus, $\tilde \tau_S$ is an $S$-optimal stopping time.

(ii) The right continuity of $(Y_t)$ and $(\xi_t)$ ensures that $Y_{\tau^*_{S}} = \xi_{\tau^*_{S}}$ a.s. 
By definition of $\tau^*_S$, we have that almost surely on $[S, \tau^*_S[$, the process $(Y_t)$ is strictly greater than the obstacle $(\xi_t)$ and hence the process $A$ is constant on $[S, \tau^*_S[$ and even on $[S, \tau^*_S]$ because $A$ is continuous (see Theorem \ref{exiuni}). We derive that $(Y_s, S\leq s \leq 
\tau^*_S)$ is the solution of the BSDE associated with terminal time $\tau^*_{S}$, terminal
condition $\xi_{\tau^*_{S}}$ and driver $f$, that is, 
$
Y_s = X_s (\xi_{\tau^*_S}, \tau^*_S)
$,  $S \leq s  \leq \tau^*_S$ a.s.\,
Hence,
$\tau_S^*$ is an  $S$-optimal stopping time.\\
Furthermore, for each $\varepsilon >0$ , $\tau^\varepsilon_S \leq \tau^*_S$ a.s. By letting $\varepsilon$  tend to $0$, we get $\tilde \tau_S \leq \tau^*_S$. a.s.

(iii) Let $\hat{\tau}$ be an $S$-optimal stopping time. By the strict comparison theorem for non reflected BSDEs (or Proposition 
\ref{optcri}), we have 
$Y_{\hat{\tau}} = \xi_{\hat{\tau}} $ a.s.
Hence, by definition of $\tau_S^*$, we have
$ \hat{\tau} \geq \tau_S^*$ a.s.\, Thus,  $\tilde \tau_S \geq \tau^*_S$ a.s.\,, which, with the other inequality, yields
 that $\tilde \tau_S = \tau^*_S$ a.s. We also have proven that $\tau^*_S$ is the minimal $S$-optimal stopping time.
\fproof

\begin{remark} Consider the case of a Brownian filtration and a continuous obstacle $(\xi_t)$. The second assertion of the above theorem concerning the optimality of $\tau^*_S$, corresponds to Theorem 5.9 in El Karoui and Quenez (1996).
\end{remark}

\section{Comparison theorems for RBSDEs with jumps and optimization problems}\label{sec-charc} 
\subsection{Comparison theorems for RBSDEs with jumps}
We  now state a comparison theorem for RBSDEs with jumps.
\begin{Theorem}[Comparison theorem for RBSDEs.]\label{thmcomprbsde}
\label{comparison} Let $\xi^{1}$, $\xi^{2}$ be two obstacle processes in ${\cal S}^2$. Let $f^1$and $f^2$ 
 be  Lipschitz drivers. Suppose that $f^1$   satisfies Assumption~\eqref{Royer}.
Suppose  that
\begin{itemize}
\item   $\xi_t ^{2}\le \xi_t ^{1}$, $0\leq t \leq T$ a.s.

\item  $f^{2}(t,y,z,k) \le f^{1}(t,y,z,k),\,\, \text{ for all  }  (y,z,k) \in \R^2 \times {\cal L}^2_\nu; \;\;  dP\otimes dt-a.s.\,$
\end{itemize}
Let $(Y^{i}, Z^{i}, k^i, A^i)$  be
the solution of the RBSDE associated with $(\xi ^{i},f^{i})$ , $i=1,2$. Then, 
$$Y_{t}^{2}\le
Y_{t}^{1}, \,\,\forall t\in [0,T] \text{ a.s. } $$
\end{Theorem}

\dproof We give here a  simple proof 
based on the characterization of solutions of RBSDEs (Theorem \ref{caracterisation}) and 
on the comparison theorem for non reflected BSDEs.
 Let $t$ $\in$ $[0,T]$. For each $\tau \in T_{t}$, let us denote by $X^{i}( \xi^i 
_{\tau}, \tau)$ the unique solution 
 of 
the BSDE associated with $(\tau, \xi^i _{\tau}, f^{i})$ for $i=1,2$.
By the comparison theorem  for BSDEs with jumps \cite{QuenSul}, the following inequality 
$$X_t^{2}(\xi^2 _{\tau}, \tau) \leq X_t^{1}( \xi^1 
_{\tau}, \tau) \,\,a.s.\,$$
holds for each $\tau$ in $\T_t$.
Hence, by taking the essential supremum over $\tau$ in $\T_t$ and  using Theorem \ref{caracterisation}, we get
$$Y^{2}_{t}= {\rm ess} \sup_{\tau \in \T_t}X^2_t(\xi^2_{\tau}, \tau) \leq {\rm ess} \sup_{\tau \in \T_t}X^1_t(\xi^1_{\tau}, \tau)=Y^{1}_{t}\,\,a.s.  $$
\fproof

\begin{remark}
The result still holds when $f^2$ satisfies Assumption~\eqref{Royer} instead of $f^1$.
\end{remark}

 We now provide a strict comparison theorem. The first assertion addresses the particular case when the obstacle is left-upper semicontinuous along stopping times and the second one deals with the general case.
\begin{theorem}[Strict comparison.]\label{sctun} 
Suppose that  the assumptions of the comparison theorem (Th. \ref{thmcomprbsde}) hold and that the driver $f^1$
satisfies Assumption \ref{Royer} with 
\begin{equation}\label{hypsupbis}
\theta_t^{x,\pi,l_1,l_2} > -1 \quad dt \otimes dP-{\rm a.s.}
\end{equation}
 Let $S$ in $\T_0$ and suppose that 
$Y^1_{S} = Y^2_{S}$ a.s.
\begin{enumerate}
\item
Suppose that  $\xi^1$ and $\xi^2$ are left-upper semicontinuous along stopping times.
Let $\tau^*_i = \tau^*_{i,S}:= \inf \{s \geq S ; \, Y^i_s = \xi^i_s \}, \; i = 1,2.$ Then, 
$$Y^1_{t} =Y^2_{t} ,\;\;  S \leq t \leq \tau_1^* \wedge \tau_2^* \; \text{ a.s}. $$
and
\begin{equation}\label{autre}
f^2(t,Y^2_t,Z^2_t,k^2_t) \, =\,  f^{1}(t,Y^2_t,Z^2_t,k^2_t) \, \, \,\,\,  \;\;   S \leq t \leq  \tau_1^* \wedge \tau_2^*, \,\,
dP\otimes dt-a.s.
\end{equation}
Moreover if $\xi^1 = \xi^2$ a.s., then $\tau_1^* = \tau_2^*$ a.s. and
 $ Y^1_{\tau_1^*} = Y^2_{\tau_1^*} = \xi^1_{\tau_1^* } $ a.s
 \item 
Consider the general case where $\xi^1$ and $\xi^2$ are not supposed to be left-upper semicontinuous along stopping times.
For  $\varepsilon >0$, define
$$\tau^{\varepsilon}_i:= \inf \{ t \geq S,\,\, Y^i_t \leq \xi^i_t + \varepsilon \} \;\; 
\text{ and  } \;\; \tilde \tau_i:= \lim_{\varepsilon \downarrow 0} \uparrow \tau^{\varepsilon}_i \;\; i = 1,2.$$ 
Then, for each $\varepsilon >0$, 

\begin{equation}\label{unoss}
Y^1_{t} =Y^2_{t} , \;\;  S \leq t \leq  \tau^{\varepsilon}_1 \wedge \tau^{\varepsilon}_2. \quad \text{a.s.}
\end{equation}
 Moreover,
$$
f^2(t,Y^2_t,Z^2_t,k^2_t) \, =\,  f^{1}(t,Y^2_t,Z^2_t,k^2_t) \, \, \,\,\,  \;\;   S \leq t \leq  \tilde \tau_1 \wedge \tilde \tau_2  , \,\,
dP\otimes dt-a.s.
$$
and if $\xi^1 = \xi^2$ a.s., then for each $\varepsilon >0$, 
$\tau^{\varepsilon}_1 = \tau^{\varepsilon}_2$ a.s. and 
$\tilde \tau_1 = \hat \tau_2 $.
\end{enumerate}
\end{theorem}

\dproof Suppose that  $\xi^1$ and $\xi^2$ are left-upper semicontinuous along stopping times.
 By the existence theorem (see Theorem \ref{existence}), $\tau^*_1$ is optimal for Problem~(\ref{prixam}) with $f= f^1$, $\xi=\xi^1$, that is 
$$Y^1_S =  {\rm ess } \sup_{\tau \in \T_S} X^1_S(\xi^1_\tau, \tau) = 
X^1_S( \xi^1_{\tau^*_1},\tau^*_1) \quad \text{a.s}$$
where $X^1( \xi^1_{\tau^*_1},\tau^*_1)$ denotes the solution of the BSDE associated with terminal time $\tau^*_1$, terminal
condition $\xi^1_{\tau^*_1}$ and driver $f^1$.
Furthermore,  by Theorem \ref{existence}, we have $$Y^1_t = X^1_t( \xi^1_{\tau^*_1},\tau^*_1), \,\,\, S \leq t \leq \tau^*_1 \,\,\,{\rm a.s.}$$ 
Moreover $\tau^*_2$ is optimal for Problem~(\ref{prixam}) with $f= f^2$, $\xi=\xi^2$,  that is, 
$$Y^2_S = {\rm ess}  \sup_{\tau \in \T_S} X^2_S(\xi^2_\tau, \tau) = 
X^2_S( \xi^2_{\tau^*_2},\tau^*_2),$$
where $X^2( \xi^2_{\tau^*_2},\tau^*_2)$ denotes the solution of the BSDE associated with terminal time $\tau^*_2$, terminal
condition $\xi^2_{\tau^*_2}$ and driver $f^2$.\\
Also, $Y^2_t = 
X^2_t( \xi^2_{\tau^*_2},\tau^*_2)$ $S \leq t \leq \tau^*_2$ a.s.
Hence $$Y^1_t = X^1_t( Y^1_{{\tau^*_1}\wedge \tau_2^*},\tau^*_1\wedge \tau_2^*), \;
\text{and} \; 
Y^2_t = X^2_t( Y^2_{{\tau^*_1} \wedge {\tau_2^*}},\tau^*_1 \wedge \tau_2^*),\,\,\,S \leq t \leq \tau^*_1 \wedge \tau_2^*\,\,\,{\rm a.s.}$$

Since  $f^1 \geq f^2$ and  $\xi^1 \geq \xi^2$,  the comparison theorem for RBSDEs (Th. \ref{thmcomprbsde}) yields 
that  $Y^1_{\tau^*_1 \wedge \tau_2^*}\geq Y^2_{\tau^*_1 \wedge \tau_2^*}$  a.s.\,
By hypothesis, $Y^1_S = Y^2_S$. Now, Assumption (\ref{hypsupbis}) allows us to apply the strict comparison theorem for non reflected BSDEs with jumps (see \cite{QuenSul} Th 4.4) for terminal time $\tau^*_1 \wedge \tau_2^*$. Hence, we get 
$Y^1_{t} =Y^2_{t} ,\;\;  S \leq t \leq \tau_1^* \wedge \tau_2^* \;$ a.s.\,, and equality (\ref{autre}), which provides the desired result.

Suppose now that  $\xi^1 = \xi^2 = \xi$ a.s. Then, using  $Y^2 \leq Y^1$, we get 
$\tau^*_2 \leq \tau^*_1$ a.s.
Since we have already shown that $Y^1_{\tau^*_2} = Y^2_{\tau^*_2}$ a.s., 
and since $   Y^2_{\tau^*_2} = \xi_{\tau^*_2} $ a.s.
Hence $Y^1_{\tau^*_2} = \xi_{\tau^*_2}$ and  $\tau^*_1 \leq  \tau^*_2$ a.s. 
It follows that $\tau^*_1 =  \tau^*_2$ a.s.

Let us now consider the general case where the obstacles are not supposed to be left-upper semicontinuous along stopping times. \\
Let $\varepsilon $ $>0$.
By a property of $\tau^{\varepsilon}_1$ (see Lemma \ref{lala}), we have 
$$Y^1_t = X^1_t( Y^1_{\tau^{\varepsilon}_1},\tau^{\varepsilon}_1), \,\,\,S \leq t \leq \tau^{\varepsilon}_1 \,\,\,{\rm a.s.}\,$$
Similarly,
 $Y^2_t = 
X^2_t( Y^2_{\tau^{\varepsilon}_2},\tau^{\varepsilon}_2), S \leq t \leq \tau^{\varepsilon}_2$ a.s.
By the same arguments as above with $\tau^*_1$ and $\tau^*_2$ replaced by $\tau^{\varepsilon}_1$ and 
$\tau^{\varepsilon}_2$ respectively, we derive the desired result.

Suppose now that  $\xi^1 = \xi^2 = \xi$ a.s. Since 
 $Y^2 \leq Y^1$, we have $\tau^{\varepsilon}_2 \leq \tau^{\varepsilon}_1$ a.s.
Moreover by inequality Lemma \ref{lala} and Assumption (\ref{unoss}), we have 
$$\xi_{\tau_2^{\varepsilon}} + \varepsilon  \geq Y^2_{\tau_2^{\varepsilon}} = 
 Y^1_{\tau_2^{\varepsilon}}  \quad \mbox{a.s.}$$
Consequently, 
$\tau_2^{\varepsilon} \geq \tau_1^{\varepsilon}$ a.s.\, Since we have already shown the converse inequality, we have $\tau_2^{\varepsilon} = \tau_1^{\varepsilon}$ a.s. 

\fproof

\subsection{Optimization problems for RBSDEs }\label{mixed1}

We use the following setup: 
Let $\xi$ in ${\cal S}^2$ and let $(f, f^{\alpha}; \alpha \in \Ac)$ be a family of Lipschitz drivers satisfying  Assumption \eqref{Royer}. In \eqref{Royer}, the coefficient associated with $f^{\alpha}$ (resp. $f$), is denoted by
$\theta^{\alpha,x, \pi,l}$ (resp. $\theta^{x, \pi,l}$).

We denote by  $(Y,Z,k)$  the solution of the RBSDE associated to obstacle $(\xi_t)$ and driver $f$, and by  $(Y^{\alpha},Z^{\alpha}, k^{\alpha})$  the solution of the RBSDE associated with obstacle $(\xi_t)$ and driver $f^{\alpha}$.
Also, for each $\tau  \in \T_0$ and $\zeta \in L^2(\FC_{\tau})$, we denote by 
$(X(\zeta, \tau), \pi(\zeta, \tau), l(\zeta, \tau))$ the solution of the BSDE  associated with driver $f$, terminal conditions $\zeta$, $\tau$, and by 
  $(X^\alpha(\zeta, \tau), \pi^\alpha (\zeta, \tau), l^\alpha (\zeta, \tau))$ the solution of the BSDE  associated with driver $f^\alpha$ and  terminal conditions $\zeta$, $\tau$. 

 From the comparison theorem, we derive a first optimization principle for RBSDEs which  generalizes   the result established by El Karoui and Quenez in \cite{EQ96} to the case of jumps. 
 \begin{Proposition} [Optimization principle for RBSDEs I] \label{minima}
Suppose that 
\begin{enumerate}
\item  For each $\alpha$ $\in$ $\Ac$, $f(t,y,z,k) \le f^{\alpha}(t,y,z,k),\,\, \text{ for all  }  (y,z,k) \in   \R^2 \times {\cal L}^2_\nu; \;\;  dt \otimes dP-a.s.\,$
\item
There exists $\bar \alpha$ $\in$ $\Ac$ such that
\begin{equation}\label{existalp}
f(t,Y_t, Z_t,k_t) = {\rm ess} \inf_{\alpha}f^{\alpha}(t,Y_t, Z_t,k_t) = 
f^{\bar \alpha}(t,Y_t, Z_t,k_t),  \; \; 0 \leq t \leq T, \;\;   dt \otimes dP-{\rm a.s.}
\end{equation}
\end{enumerate}
Then, for each $S \in \T_0$,
\begin{eqnarray}\label{cequal}
Y_S = {\rm ess} \inf_{\alpha}Y^{\alpha}_S = Y^{\bar \alpha}_S \quad {\rm a.s.}
\end{eqnarray}
\end{Proposition}

\dproof For each $\alpha$, since Condition 1.\, is satisfied and, since $f^{\alpha}$ satisfies Assumption \ref{Royer}, the comparison theorem for RBSDEs yields  (see Theorem \ref{thmcomprbsde}) that  
$Y \leq Y^{\alpha}$. It follows 
that  for each $S \in \T_0$,
$$
Y_S \leq {\rm ess} \inf_{\alpha}Y^{\alpha}_S  \quad {\rm a.s.}
$$
Now, by condition 2.\,, $Y $ is a solution of the RBSDE associated with $f^{\bar \alpha}$. By uniqueness of the solution of this RBSDE, we have 
 $Y = Y^{\bar 
\alpha}$,  which leads to 
equality \eqref{cequal}. 
\fproof

\begin{remark}
This Proposition still holds if $f$ does not satisfy Assumption \eqref{Royer}.
\end{remark}

  \begin{Proposition}[Optimization principle for RBSDEs II]\label{2principle}
Suppose that the drivers $f^{\alpha}$, $\alpha$ $\in$ $\Ac$ satisfy $f \le f^{\alpha}$ and are equi-Lipschitz with constant $C$.\\
Suppose moreover that for each $\eta >0$ , there exists $\alpha^{\eta}$ $\in$ $\Ac$ such that
\begin{equation}\label{ep1}
 f(t,Y_t,Z_t,k_t) \geq f^{{\alpha^{\eta}}}(t,Y_t,Z_t,k_t)- \eta,  \, \;\   0 \leq t \leq  T, \,
 dP\otimes dt-{\rm a.s.}
\end{equation}
Then, for each $S$ $\in$ $\T_0$, we have
\begin{equation}\label{q}
Y_S={\rm ess} \inf_{\alpha} Y^{\alpha}_S \quad {\rm a.s.}
\end{equation}
\end{Proposition}

\dproof
 Since $f \le f^{ \alpha}$, we have $Y \leq Y^{ \alpha}$ a.s.\,for each $\alpha \in \Ac$. It follows that for each 
 $S$ $\in$ $\T_0$, we have 
 $Y_S\leq {\rm ess} \inf_{\alpha} Y^{\alpha}_S $ a.s.\, 
Since Assumption (\ref{ep1}) holds, by using estimation \eqref{A26}, with $\eta = \frac{1}{C^2}$ and $\beta = 3C^2 + 2C$, we derive 
that there exists a constant $K\geq 0$, which depends only on $C$ and $T$, such that, for each $\eta>0$  and for each $S$ $\in$ $\T_0$,
$$Y_S + K \,\eta \geq Y^{{\alpha^{\eta}}}_S \geq {\rm ess} \inf_{\alpha} Y^{\alpha}_S \quad {\rm a.s.}
$$
Equality (\ref{q}) thus follows.
\fproof

By using the strict comparison theorem for reflected BSDEs (see Theorem \ref{sctun}), we provide some necessary and sufficient conditions of optimality at a given time $S$ $\in$ $\T_0$. 
 
 \begin{Theorem}[Optimality criteria for RBSDEs.]\label{firstcrit}
Suppose that  for each $\alpha$ $\in$ $\Ac$, $f \leq f^{\alpha}$.
Let ${\bar \alpha} \in \Ac$, and suppose that 
 in Assumption \ref{Royer}  the coefficient $\theta_{\bar \alpha}$ corresponding to driver $f^{\bar \alpha}$ satisfies
$\theta^{\bar \alpha,x, \pi,l} > -1$, for each $x, \pi,l$.
\begin{enumerate}
\item
Suppose that the obstacle $\xi$ is left-upper semicontinuous along stopping times.
Define for each $S$ in $\T_0$, 
$$\tau^*_S:= \inf \{ t \geq S,\,\,  Y_t = \xi_t\}.$$
The  parameter  $\bar \alpha$ is $S$-optimal (i.e. 
$ess \inf_{\alpha} Y^{\alpha}_S =  Y^{{\bar \alpha}}_S $) if and 
only if 
\begin{eqnarray}\label{strict}
Y^{\bar \alpha}_{ \tau^*_S}= \xi_{ \tau^*_S}\, \, {\rm a.s.} \,; \;  f(t,Y_t,Z_t,k_t) = f^{{\bar \alpha}}(t,Y_t,Z_t,k_t),  \;   S \leq t \leq   \tau^*_S, \,
dP\otimes dt-a.s.
\end{eqnarray}
\item Consider the general case when the obstacle is not supposed to be left-upper semicontinuous along stopping times. Define  for each $\varepsilon >0$, and each $S$ $\in$ $\T_0$,  the stopping time $\tau^{\varepsilon}_S:= \inf \{ t \geq S,\,\, Y_t \leq \xi_t + \varepsilon\}$.\\
A parameter  $\bar \alpha$ is $S$-optimal (i.e. 
$ess \inf_{\alpha} Y^{\alpha}_S =  Y^{{\bar \alpha}}_S $) if and only if
for each 
$\varepsilon >0$, 
\begin{eqnarray}\label{strictbis}
Y^{\bar \alpha}_{ \tau^{\varepsilon}_S}\leq  {\xi_{\tau^{\varepsilon}_S}} + {\varepsilon}\, \, {\rm a.s.} \,; \;  f(t,Y_t,Z_t,k_t)  = f^{{\bar \alpha}}(t,Y_t,Z_t,k_t),  \,   S \leq t \leq  \tau^{\varepsilon}_S, \,
 dP\otimes dt-{\rm a.s.}
\end{eqnarray}
\end{enumerate}
Also, in both cases, $Y_S=ess \inf_{\alpha} Y^{\alpha}_S =  Y^{{\bar \alpha}}_S $ a.s.
\end{Theorem}

\begin{remark}\label{rem63}
Note that  in  the first assertion,  even if the assumption $\theta^{\bar \alpha,x, \pi,l} > -1$ is not satisfied, \eqref{strict} implies that 
$\bar \alpha$ is $S$-optimal. 
The same holds for assertion 2.
\end{remark}

\dproof {\em 1.}\, Suppose that $\bar \alpha$ is $S$-optimal. Note that, since $Y \leq Y^{{\bar \alpha}}$, it follows that 
$\tau^*_S \leq  \tau^{\bar \alpha, *}_S$ where 
$\tau^{\bar \alpha, *}_S :=  \inf \{ t \geq S,\,\,  Y^{{\bar \alpha}}_t = \xi_t\}$.
By the first strict comparison theorem for RBSDEs (Theorem \ref{sctun} 1.) applied to $\xi^1= \xi^2 =\xi$, $f^1=f$, $f^2= f^{{\bar \alpha}}$, $Y^1=Y$, $Y^2= Y^{{\bar \alpha}}$, we derive that equalities (\ref{strict}) hold.
 
It remains to show the converse. Suppose that equalities (\ref{strict}) hold. Then, by the optimality of $\tau^*_S$ for $Y_S$, we have
$$Y_t=   X_t(\xi_{\tau^*_S}, \tau^*_S), \, \, \,\,\,  \;\;   S \leq t \leq   \tau^*_S, \,\, {\rm a.s.}$$
This with equality (\ref{strict}) and the uniqueness result for BSDEs
 leads to 
$$Y_t= X_t(\xi_{\tau^*_S}, \tau^*_S)=X^{\bar \alpha}_t(\xi_{\tau^*_S}, \tau^*_S)
=X^{\bar \alpha}_t(Y^{\bar \alpha}_{ \tau^*_S}, \tau^*_S), \, \, \,\,\,  \;\;   S \leq t \leq   \tau^*_S, \,\, {\rm a.s.}\,,$$
Moreover, according to the previous equalities, $X^{\bar \alpha}_t(Y^{\bar \alpha}_{ \tau^*_S}, \tau^*_S)
=Y_t \geq \xi_t$, $\,S \leq t \leq   \tau^*_S$ a.s.\, By the 
uniqueness result for RBSDEs, it follows that 
$$Y_t =X^{\bar \alpha}_t(Y^{\bar \alpha}_{ \tau^*_S}, \tau^*_S)= Y_t ^{\bar \alpha}, \, \, \,\,\,  \;\;   S \leq t \leq   \tau^*_S, \,\, {\rm a.s.}\,$$
 By taking $t=S$, we get  $Y_S= {\rm ess}\inf_{\alpha} Y^{\alpha}_S =  Y^{{\bar \alpha}}_S$ a.s.\,, which ends the proof of the first assertion.
 
{\em 2.}\,  Suppose that $\bar \alpha$ is $S$-optimal. 
Let $$\tau^{\bar \alpha,  \varepsilon}_S :=  \inf \{ t \geq S,\,\,  Y^{{\bar \alpha}}_t \leq \xi_t + \varepsilon\}.$$
Since $Y \leq Y^{{\bar \alpha}}$, it follows that for each $\varepsilon >0$, we have 
$\tau^ \varepsilon_S \leq  \tau^{\bar \alpha, \varepsilon}_S$ a.s.\,
By the second strict comparison theorem for RBSDEs  (Theorem \ref{sctun} 2.) applied to $\xi^1= \xi^2 =\xi$, $f^1=f$, $f^2= f^{{\bar \alpha}}$, $Y^1=Y$, $Y^2= Y^{{\bar \alpha}}$, we derive that $Y^{\bar \alpha}_{\tau^{\varepsilon}_S}= Y_{\tau^{\varepsilon}_S} \leq {\xi_{\tau^{\varepsilon}_S}} + {\varepsilon}$ a.s.\,and  $f(t,Y_t,Z_t,k_t)  = f^{{\bar \alpha}}(t,Y_t,Z_t,k_t)$, $ S \leq t \leq  \tau^{\varepsilon}_S$, $dP\otimes dt$-a.s.

It remains to show the converse. Suppose that equalities (\ref{strictbis}) hold. Note first that since $f \le f^{\bar \alpha}$, we clearly have $Y_S \leq Y^{\bar \alpha}_S$ a.s.\\
Let us now show  that $Y_S \geq Y^{\bar \alpha}_S$ a.s.\, By a property of $\tau^{\varepsilon}_S$ 
(see Lemma \ref{lala}), we have 
$$Y_t = X_t(Y_{\tau^{\varepsilon}_S}, \tau^{\varepsilon}_S), \, \, \,\,\,  \;\;   S \leq t \leq   \tau^{\varepsilon}_S, \,\, {\rm a.s.}\,,$$
Hence, using equality (\ref{strictbis}), we derive that
$$Y_t=X_t(Y_{\tau^{\varepsilon}_S}, \tau^{\varepsilon}_S)= X^{\bar \alpha}_t(Y_{\tau^{\varepsilon}_S}, \tau^{\varepsilon}_S), \, \, \,\,\,  \;\;   S \leq t \leq    \tau^{\varepsilon}_S, \,\, {\rm a.s.}\,.$$
By the comparison theorem for non reflected BSDEs and the inequality 
$Y_{\tau^{\varepsilon}_S}\geq \xi_{\tau^{\varepsilon}_S}$ a.s.\,, we have 
$$Y_t =X^{\bar \alpha}_t(Y_{\tau^{\varepsilon}_S}, \tau^{\varepsilon}_S) \geq  X^{\bar \alpha}_t(\xi_{\tau^{\varepsilon}_S}, \tau^{\varepsilon}_S), \, \, \,\,\,  \;\;   S \leq t \leq    \tau^{\varepsilon}_S,\quad {\rm a.s.}$$

Now, by the a priori estimates (see \cite{QuenSul}), we have
$$Y_S \geq X^{\bar \alpha}_S(\xi_{\tau^{\varepsilon}_S}, \tau^{\varepsilon}_S)\geq X^{\bar \alpha}_S(\xi_{\tau^{\varepsilon}_S} + \varepsilon\,, \tau^{\varepsilon}_S) - \varepsilon e ^{ \frac{\beta T } {2}}\quad {\rm a.s.}\,$$
with $\beta = 3 C^2 + 2C$, where  $C$ is the Lipschistz constant of $f^{\bar \alpha}$.
Since by assumption, $  \xi_{\tau^{\varepsilon}_S}+\varepsilon \geq Y^{\bar \alpha}_{ \tau^{\varepsilon}_S}\,$ a.s.\,, the comparison theorem for non reflected BSDEs yields that

$$  Y_S + \varepsilon e ^{ \frac{\beta T } {2}} \geq  X^{\bar \alpha}_S(\xi_{\tau^{\varepsilon}_S} + \varepsilon, \tau^{\varepsilon}_S) 
\geq X^{\bar \alpha}_S(Y^{\bar \alpha}_{ \tau^{\varepsilon}_S}\,, \tau^{\varepsilon}_S) \quad {\rm a.s.}$$
Since $Y_{.}\leq Y_{.}^{\bar \alpha}$, we have $\tau^{\varepsilon}_S \leq \tau^{\bar \alpha,\varepsilon}_S$ a.s.\, (actually equality holds). Now, by Lemma \ref{lala}, the non decreasing process associated with 
$Y_{.}^{\bar \alpha}$ is constant on $[S, \tau^{\bar \alpha ,\varepsilon}_S]$ and hence on 
$[S, \tau^{\varepsilon}_S]$. Thus, $(Y^{\bar \alpha}_t,  S \leq t \leq \tau_S^\varepsilon)$ is the solution of the non reflected BSDE associated 
with driver $f^{\bar \alpha}$, terminal time $ \tau_S^\varepsilon$, and terminal condition $Y_{ \tau_S^\varepsilon}^{\bar \alpha}$. We thus get
$$X^{\bar \alpha}_S(Y^{\bar \alpha}_{ \tau^{\varepsilon}_S}, \tau^{\varepsilon}_S)=Y^{\bar \alpha}_S\quad {\rm a.s.}$$
Consequently, for each $\varepsilon >0$, we have $  Y_S + \varepsilon e ^{ \frac{\beta T } {2}} \geq Y^{\bar \alpha}_S$ a.s.\,, and hence, $Y_S \geq Y^{\bar \alpha}_S$ a.s.\, We thus have $Y_S = Y^{\bar \alpha}_S$ a.s.\,, which provides the desired result.
\fproof

 \section{Robust optimal stopping problem}\label{mixed}

We  now consider the optimal stopping problem when there is ambiguity on the risk-measure modeling. Let $\{f^{\alpha}, \alpha \in {\cal A}\}$ be  a given 
family of Lipschitz drivers satisfying Assumption \eqref{Royer}.  
For each $\alpha$ $\in$ ${\cal A}$, let $\rho^{\alpha}$ be the risk measure induced by the BSDE with driver $f^{\alpha}$, defined as follows: for each terminal time $\tau  \in \T_0$ and position $\zeta \in L^2(\FC_{\tau})$, set 
  \begin{equation*}
  \rho^{\alpha}_t (\zeta, \tau) :=  -X^{\alpha}_t(\zeta,\tau), \,\, \,\,\,0\leq t \leq T,
  \end{equation*}
 where $X^{\alpha}_t(\zeta,\tau)$ denotes the solution of the BSDE associated with driver $f^{\alpha}$, terminal condition $\zeta$ and terminal time  $\tau$. We consider an  agent who is averse 
to ambiguity, and we define her risk measure of position $\zeta$, at each time $S$ in $\T_0$ with $S \leq \tau$ a.s.\,, as the  supremum over $\alpha$ of the associated risk-measures $\rho^{\alpha}_S (\zeta, \tau)$ that is, 
   \begin{equation*}\label{definitionbisbis}
 \text{ess} \sup_{\alpha \in  {\cal A}  } \rho^{\alpha}_S (\zeta, \tau) =  \text{ess} \sup_{\alpha \in  {\cal A}} -X^{\alpha}_S(\zeta,\tau).
  \end{equation*}
 Let $(\xi_t)$ be a dynamic position, given by an RCLL adapted process $(\xi_t)$ in ${\cal S}^2$.
 At time $S$ $\in$ $\T_0$, the agent wants to choose a stopping time $\tau$ $\in$ $\T_S$ 
which minimizes her risk measure. At time $S$, her value function is defined as
 \begin{equation}\label{2problem}
u(S):=\text{ess} \inf_{\tau \in \T_S} \text{ess}\sup_{\alpha \in  {\cal A}} \rho^{\alpha}_S (\xi_{\tau}, \tau).
\end{equation}

This leads to the following game problem.\\
Let $S \in \T_0$. Define the {\em first value function at time $S$} as
\begin{equation}\label{vdessous}
\underline V(S):={\rm ess} \inf_{\alpha \in {\cal A}}{\rm ess} \sup_{\tau \in \T_S}   X^{\alpha}_S(\xi_{\tau}, \tau) ,
\end{equation}
and the {\em second value function at time $S$} as
\begin{equation}\label{vdessus}
\bar V(S):={\rm ess} \sup_{\tau \in \T_S} {\rm ess} \inf_{\alpha \in {\cal A}}  X^{\alpha}_S(\xi_{\tau}, \tau).
\end{equation}
Note that $\bar V(S)= - u(S)$ a.s.\\
By definition, we say that there exists a {\em value function} at time $S$ for the game problem if $\bar V(S) = \underline V(S)$ a.s.\,\\
We introduce the definition of an $S$-saddle point: 
\begin{Definition}\label{defsaddle}
Let $S$ $\in$ $\T_0$. A pair  $ (\hat \tau, \hat \alpha)$ 
$\in \T_S \times {\cal A}$ is called a  {\em $S$-saddle point} if 
\begin{itemize}
\item $\bar V(S) = \underline V(S)$ a.s.\,, 
\item  the essential infimum in (\ref{vdessous})  is attained at $\hat \alpha$, 
\item  the essential supremum in (\ref{vdessus}) is attained at $\hat \tau$.
\end{itemize}
By classical results, for each $S \in \T_0$,
 $ (\hat \tau, \hat \alpha)$ is a $S$-saddle point if and only if
for each $ (\tau, \alpha)$   $\in \T_S \times {\cal A}$, 
\begin{equation}\label{classiquepointselle}
X^{\hat \alpha}_S(\xi_{\tau}, \tau) \,  \leq  \,\,\, X^{\hat \alpha}_S(\xi_{\hat \tau}, \hat \tau) \leq \,\,\, X^{ \alpha}_S(\xi_{\hat \tau}, \hat \tau) \,\,\, {\rm a.s.}
\end{equation}
\end{Definition}

Note that for each $S$ $\in$ $\T_0$, the inequality $\bar V(S) \leq \underline V(S)$ a.s.\,clearly holds.
We want to determine when the equality holds, characterize the value function, and  address the question of existence of a $S$-saddle point. 

\begin{remark}\label{SP}
If $ (\hat \tau, \hat \alpha)$ is  an $S$-saddle point, then $\hat \tau$ and $\hat \alpha$ attain respectively the 
infimum and the supremum in $\bar V (S)$ that is, 
$$\bar V (S) = {\rm ess} \sup_{\tau \in \T_S} {\rm ess} \inf_{\alpha}  X^{\alpha}_S(\xi_{\tau}, \tau) =  {\rm ess} \inf_{\alpha}  X^{\alpha}_S(\xi_{\hat \tau}, \hat \tau) = X^{\hat \alpha}_S(\xi_{\hat \tau}, \hat  \tau). $$
Hence, $\hat \tau$ is an optimal stopping time for 
the agent who wants to minimize over stopping times her 
risk-measure at time $S$ in the case of ambiguity (see \eqref{2problem}). 
Also, since $\hat \alpha$ attains the essential infimum in (\ref{vdessous}), $\rho^{\hat \alpha}$ can be interpreted as the``worst" risk measure.
\end{remark}

We will now relate the game problem to an optimization problem for RBSDEs.\\ 
Let $(Y^{\alpha},Z^{\alpha}, k^{\alpha})$ be the solution of the RBSDE with obstacle $(\xi_t)$ and driver $f^{\alpha}$.
For each $\tau  \in \T_0$ and $\zeta \in L^2(\FC_{\tau})$, let 
$(X^\alpha(\zeta, \tau), \pi^\alpha (\zeta, \tau), l^\alpha (\zeta, \tau))$ be the solution of the BSDE  with driver $f^{\alpha}$ and terminal conditions $(\zeta, \tau)$. 
  
By the characterization of RBSDEs (see Theorem~\ref{caracterisation}), for each $S \in \T_0$, we have\\
$Y_S^{\alpha}= {\rm ess} \sup_{\tau \in \T_S}   X^{\alpha}_S(\xi_{\tau}, \tau) $ a.s.\,It follows that 
\begin{equation}
\underline V(S) = {\rm ess} \inf_{\alpha \in {\cal A}} Y_S^\alpha  \quad \text{a .s. }
 \end{equation}

By using the previous results on RBSDEs, we provide the following theorem, which holds for a general adapted RCLL obstacle process $(\xi_t)$.

Let $f$ be a Lipschitz driver satisfying Assumption \eqref{Royer}. Let $(Y,Z,k)$  be the solution of the RBSDE with obstacle $(\xi_t)$ and driver $f$.
For each $\tau  \in \T_0$ and $\zeta \in L^2(\FC_{\tau})$, let 
$(X(\zeta, \tau), \pi(\zeta, \tau), l(\zeta, \tau))$ be the solution of the BSDE  with driver $f$ and terminal conditions $(\zeta, \tau)$. 

\begin{Theorem}[Verification theorem I]\label{general}
Suppose that the drivers $f^{\alpha}$, $\alpha$ $\in$ $\Ac$ satisfy $f \le f^{\alpha}$ and are equi-Lipschitz with constant $C$. Suppose that 
there exists  $\bar \alpha$ such that 
\begin{equation}\label{existalpter}
f(t,Y_t, Z_t,k_t) = {\rm ess} \inf_{\alpha \in \Ac} f^{\alpha}(t,Y_t, Z_t,k_t) = 
f^{\bar \alpha}(t,Y_t, Z_t,k_t),  0 \leq t \leq T, \;\;   dt \otimes dP-a.s.
\end{equation}
Then, there exists a value function, which is characterized as the solution of the RBSDE with obstacle $(\xi_t)$ and driver $f$, that is,
for each
$S$ $\in$ $\T_0$, we have
$$Y_S=  \underline V(S)= \bar V(S)\;\; {\rm a.s.}$$ 
\end{Theorem}
This theorem can be seen as a {\em verification theorem} in the following sense: 
if we are given a driver $f$ satisfying some appropriate conditions, the solution of the RBSDE with driver $f$ coincides with the value function of the game problem.

\dproof  Let $S$ $\in$ $\T_0$.
Let us  prove that $ \underline V(S) \leq  \bar V(S) $ a.s.\,
By  assumption \eqref{existalpter} and the optimization principle for RBSDEs (see Theorem~\ref{minima}), we have: 
\begin{equation}\label{equs}
\underline V(S) = {\rm ess} \inf_{\alpha \in {\cal A}} Y_S^\alpha = Y_S^{\bar \alpha} 
= Y_S \quad \text{a.s. }
 \end{equation}

Let $\varepsilon$ $>0$. 
By a property of $\tau^{\varepsilon}_S$ (see Lemma \ref{lala}), we have
$$Y_t =  X_t ( Y_{ \tau^{\varepsilon}_S },  \tau^{\varepsilon}_S ), \,S \leq t \leq \tau^{\varepsilon}_S,  \quad \mbox{a.s.}\,$$
If $(X_t, \pi_t,l_t)$ denotes the solution of the BSDE associated with driver $f$ and terminal conditions $(Y_{ \tau^{\varepsilon}_S}, \tau^{\varepsilon}_S)$, we thus have $(Y_t, Z_t,k_t)= (X_t, \pi_t,l_t)$ for $S \leq t \leq  \tau^{\varepsilon}_S$ a.s.\,
This with Assumption (\ref{existalpter}) ensures that
\begin{equation}\label{ff}
f(t,X_t, \pi_t,l_t) = {\rm ess} \inf_{\alpha \in \Ac}f^{\alpha}(t,X_t, \pi_t,l_t) = 
f^{\bar \alpha}(t,X_t, \pi_t,l_t),  S \leq t \leq  \tau^{\varepsilon}_S, \;\;   dt \otimes dP-a.s.
\end{equation}
Hence, the first optimization principle for non reflected BSDEs (see \cite{QuenSul}) can be applied. It follows that
\begin{equation}\label{optimum}
 X_S ( Y_{ \tau^{\varepsilon}_S },  \tau^{\varepsilon}_S ) = {\rm ess} \inf_{\alpha} X^{\alpha}_S 
( Y_{ \tau^{\varepsilon}_S},  \tau^{\varepsilon}_S ) \quad {\rm a.s.}
\end{equation}
Using the comparison theorem for non reflected BSDEs and the inequality 
$ Y_{ \tau^{\varepsilon}_S } \leq \xi_{ \tau^{\varepsilon}_S }+ \varepsilon$ a.s.\,, it follows that
 \begin{equation}\label{ys}
Y_S = {\rm ess} \inf_{\alpha} X^{\alpha}_S 
( Y_{ \tau^{\varepsilon}_S},  \tau^{\varepsilon}_S ) \leq {\rm ess} \inf_{\alpha} X^{\alpha}_S 
(\xi_{ \tau^{\varepsilon}_S }+ \varepsilon,  \tau^{\varepsilon}_S ).
\end{equation}
By the a priori estimates for non reflected BSDEs with jumps (see \cite{QuenSul}), for each $\varepsilon >0$ and for each $\alpha \in \Ac$, we have
\begin{equation*}  
 X^{ \alpha}_S ( \xi_{ \tau^{\varepsilon}_S }+ \varepsilon, \tau^{\varepsilon}_S) \leq X^{ \alpha}_S ( \xi_{\tau^{\varepsilon}_S}, \tau^{\varepsilon}_S) +\varepsilon e^{\frac{\beta T}{2}}\quad \mbox{a.s.}\,,
 \end{equation*}
with $\beta = 3C^2+2C$, where the constant $C$ is equal to the Lipschitz constant common to all the drivers $f^{\alpha}$, $\alpha \in \Ac$.
By taking the essential infimum over $\alpha$, we derive that for each $\varepsilon >0$,
 $${\rm ess} \inf_{\alpha}  X^{ \alpha}_S ( \xi_{ \tau^{\varepsilon}_S }+ \varepsilon, \tau^{\varepsilon}_S) \leq {\rm ess} \inf_{\alpha} X^{ \alpha}_S ( \xi_{\tau^{\varepsilon}_S}, \tau^{\varepsilon}_S) +\varepsilon e^{\frac{\beta T}{2}}
 \leq \bar V(S) +\varepsilon e^{\frac{\beta T}{2}}
 \quad \mbox{a.s.}\,,$$
 where the last inequality follows from the fact that
 $$\bar V(S) =  {\rm ess} \sup_{\tau \in \T_S}{\rm ess} \inf_{\alpha} X^{ \alpha}_S ( \xi_{\tau}, \tau)\quad \mbox{a.s.}$$
Using (\ref{ys}), we get $Y_S 
 \leq \bar V(S) +\varepsilon e^{\frac{\beta T}{2}}$ a.s.\,
  Since $\underline V(S)= Y_S$ a.s. (see \eqref{equs}), it follows that 
 for each $\varepsilon >0$, we have
\begin{equation*}\label{yy}
\underline V(S)=Y_S    \leq \bar V(S) +\varepsilon e^{\frac{\beta T}{2}}\quad \mbox{a.s.}\,
\end{equation*}
Hence, $\underline V(S)=Y_S    \leq \bar V(S)$ a.s.\, Since $\bar V(S) \leq \underline V(S)$ a.s.\,, it follows that  $\underline V(S)=Y_S    = \bar V(S)$ a.s.\,
The proof is thus complete.
\fproof

\begin{remark}\label{remark68}
Suppose that for each $\alpha$ in $\cal A$, $f \leq f^\alpha$ and  the drivers $f^{\alpha}$ are equi-Lipschitz. 
Let $S $ in $\T_0$. Assume there exists   $\bar \alpha$ such that
for each 
$\varepsilon >0$, 
\begin{eqnarray}\label{649}
Y^{\bar \alpha}_{ \tau^{\varepsilon}_S}\leq  {\xi_{\tau^{\varepsilon}_S}}+{\varepsilon}\,  {\rm a.s.} \,\,{\rm and }\,  \;\; f(t,Y_t,Z_t,k_t)  = f^{{\bar \alpha}}(t,Y_t,Z_t,k_t),  \,   S \leq t \leq  \tau^{\varepsilon}_S, \,
 dP\otimes dt-a.s.
\end{eqnarray}
Then, we have
 $$Y_S=\underline V(S) = \bar V(S).$$
Note that \eqref{649} is  weaker  than \eqref{existalpter}.
This result follows from  the second optimality criterium (see Theorem \ref{firstcrit} 2.) and the same arguments as above.  
\end{remark}

We stress on that the above theorem holds without making the left-upper semicontinuity hypothesis on $\xi$ along stopping times and hence, it may be that there does not exist any optimal stopping time for $Y_S =  {\rm ess} \sup_{\tau \in \T_S} X_S(\xi_{\tau}, \tau)$ and that there does not exist any $S$-saddle point. 

We now show the following verification theorem, which holds under weaker hypotheses. 
 \begin{Theorem}[Verification Theorem II] \label{optieps}
Suppose that  for each $\alpha$ $\in$ $\Ac$, 
$f \le f^{\alpha}$. 
Suppose that for each $\eta >0$, there exists $\alpha^{\eta}$ $\in$ $\Ac$ such that
\begin{equation}\label{epsilon}
 f(t,Y_t,Z_t,k_t) \geq f^{{\alpha^{\eta}}}(t,Y_t,Z_t,k_t)- \eta,  \, \;\   0 \leq t \leq  T, \,
 dP\otimes dt-a.s.
\end{equation}
Then, for each $S \in \T_{0}$, the equality $Y_S=\underline V(S) = \bar V(S)$ holds a.s.
\end{Theorem}

\dproof
By Theorem \ref{2principle}, we already know that
$Y_S={\rm ess} \inf_{\alpha} Y^{\alpha}_S = \underline V(S)$ a.s.\\
 Since $f \le f^{ \alpha^{\eta}}$, we  have $Y_S \leq Y^{ \alpha^{\eta}}_S$ a.s. 

For each $\varepsilon >0$, by a property of $\tau^{\varepsilon}_S$ (see Lemma \ref{lala}),
we have 
$$(Y_t, Z_t, k_t) = (X_t, \pi_t, l_t)\, \, \,\,\,  \;\;   S \leq t \leq   \tau^{\varepsilon}_S, \,\, {\rm a.s.}\,.$$
By assumption \eqref{epsilon}, we have
\begin{equation}\label{epsilonb}
 f(t,X_t,\pi_t,l_t) \geq f^{{\alpha^{\eta}}}(t,X_t,\pi_t,\eta_t)- \eta,  \, \;\   S \leq t \leq  \tau^{\varepsilon}_S, \,
 dP\otimes dt-a.s.
\end{equation}
and this holds for each $\eta \geq 0$.
By the second optimization principle for non reflected BSDE (see \cite{QuenSul}, Theorem 4.6), we have
\begin{equation*}\label{una}
Y_S = X_S = X_S(Y_{\tau^{\varepsilon}_S}, \tau^{\varepsilon}_S)= 
{\rm ess} \inf_\alpha  X^\alpha_S(Y_{\tau^{\varepsilon}_S}, \tau^{\varepsilon}_S)
 \,\, {\rm a.s.}\,
 \end{equation*}
 The end of the proof is the same as that of Theorem \ref{general}. 
 \fproof

%

From the above theorems, we derive a saddle point criterium.
\begin{Corollary} \label{scriterium} Suppose that the assumptions of Theorem \ref{general} or Theorem \ref{optieps} are satisfied. 
Let $S$ $\in$ $\T_0$. For each stopping time $\hat \tau$ $\in \T_S$ and for each $\hat \alpha$ $\in {\cal A}$, the 
pair $(\hat \tau, \hat \alpha)$ is 
an $S$-saddle point if and only if $\hat \tau$ is an optimal stopping time for 
$Y_S=  {\rm ess} \sup _{\tau \in \T_S} X_S(\xi_{\tau}, \tau)$ and $\hat \alpha$ is optimal for 
$Y_S = {\rm ess} \inf _{\alpha \in {\cal A}} Y_S^{\alpha}$.

\end{Corollary}

\dproof
By Theorem \ref{general}  or \ref{optieps}, we have 
$ \bar V(S) = \underline V(S) = Y_S $ a.s. 
The result follows from  the definition of an $S$-saddle point (see Definition \ref{defsaddle}).
\fproof

The following existence result clearly follows.

\begin{Corollary} 
\label{existenceun}
Suppose that the assumptions of Theorem \ref{general} hold and 
 that the obstacle $\xi$ is left-upper semicontinuous along stopping times. Let 
 $\tau^*_S := \inf \{u \geq S ; \, Y_u = \xi_u \}$.\\
Then, for each $S \in \T_0$,   $( \tau_S^*, \bar \alpha)$ is an $S$-saddle point. 
\end{Corollary}
\begin{remark}
 This corollary generalizes a similar result of  \cite{EQ96} obtained  in the case of a  Brownian framework and a continuous obstacle.
\end{remark}
By Theorem \ref{firstcrit} and Remark \ref{rem63}, 
we get the following existence result which holds under a weaker hypothesis. 
\begin{Corollary}
\label{existencedeux}
Suppose that the assumptions of Theorem \ref{optieps} are satisfied and that 
the obstacle $\xi$ is left-upper semicontinuous along stopping times. Let $S$ in $\T_0$.
Suppose that there exists $\bar \alpha$  such that 
\begin{equation}\label{cribis}
Y^{\bar \alpha}_{ \tau^*_S}= \xi_{ \tau^*_S}\, \, {\rm a.s.} \, \text{ and } f(t,Y_t,Z_t,k_t) = f^{{\bar \alpha}}(t,Y_t,Z_t,k_t),  \; \;  S \leq t \leq   \tau^*_S, \,\, dP\otimes dt-{\rm a.s.}
\end{equation}
Then,  $( \tau^*_S, \bar \alpha)$ is an $S$-saddle point. 
\end{Corollary}

\section{Application to the case of multiple priors}\label{application}

We now apply these results to an optimal stopping problem for dynamic risk-measures in the case of multiple priors.
Let $A$ be a Polish space (or a Borelian subset of a Polish space) and let ${\cal A}$ the set of $A$-valued predictable processes $\alpha$. With each coefficient $\alpha \in {\cal A}$, is associated a model 
via a probability measure $Q^{\alpha}$ called {\em prior} as well as a dynamic risk measure $\rho^{\alpha}$. More precisely, for each 
$\alpha \in {\cal A}$, let $Z^{\alpha}$ be the solution of the SDE: 
$$ dZ_t^{\alpha} = Z_{t^-} ^{\alpha}\left( \beta^1 (t, \alpha_t) dW_t + \int_{\R^*}
\beta^2(t, \alpha_t,u) d\tilde N(dt,du)\right)\,; \,\,\,\, Z^{\alpha}_0= 1,$$
where 
 $\beta^1: (t, \omega ,\alpha) \mapsto \beta^1(t, \omega,\alpha)$, is  a ${\cal P}\otimes {\cal B}(A)$-measurable function defined on $[0,T] \times \Omega \times  A$ and valued in $[-C,C]$, with $C >0$, 
 and  $\beta^2: (t, \omega,\alpha,u) \mapsto \beta^2(t, \omega,\alpha,u)$ is  a ${\cal P}\otimes {\cal B}(A) \otimes {\cal B}(\R^*)$-measurable function defined on $ [0,T] \times \Omega \times A \times \R^*$ which satisfies $dt \otimes dP \otimes d\nu (u)$-a.s.
\begin{equation}\label{roy}
\beta^2(t,\alpha,u) \geq C_1  \;\; \text{ and }
\;\; |\beta^2(t,\alpha,u) | \leq \psi(u),
\end{equation}
with $C_1 > -1$ and  $\psi$  is a bounded function $\in$ $L^{p}_{\nu}$ for all $p \geq 1$.
Hence, $Z_{T}^{\alpha}>0$ a.s. and, by Proposition A1 in \cite{QuenSul},  
 $Z_{T}^{\alpha}$ $\in$ $L^{p}({\cal F}_{T})$ for all $p \geq 1$. 
  
For each 
$\alpha \in {\cal A}$, let $Q^{\alpha}$  be the probability measure equivalent to 
$P$ which admits 
 $Z_{T}^{\alpha}$ as density with respect to $P$ on ${\cal F}_{T}$.
By Girsanov's theorem, the process $W^{\alpha}_t := W_t - \int_0^t \beta^1 (s, \alpha_s) ds$ is a Brownian motion under $Q^{\alpha}$ and $N$ is a Poisson random measure independant from $W^{\alpha}$ under $Q^{\alpha}$ with compensated process
$\tilde N^{\alpha}(dt,du) =\tilde N(dt,du)  -\beta^2(t,\alpha_t,u) \nu(du)dt$.

For each control $\alpha$, the associated dynamic risk measure is induced by a BSDE under $Q^{\alpha}$ and driven by $W^{\alpha}$  and $\tilde N^{\alpha}$, which makes sense since we have a $Q^{\alpha}$-martingale representation property (see Lemma 5.7 in \cite{QuenSul}).
We introduce a function\\
$F:[0,T] \times  \Omega \times \R \times L^2_\nu \times A \rightarrow \R $ ; 
$(t, \omega, \pi, \ell, \alpha) \mapsto  F(t, \omega, \pi, \ell, \alpha) $
which  is $ {\cal P} \otimes {\cal B}(\R)  \otimes {\cal B}(L^2_\nu)  \otimes {\cal B}(A)$-measurable. 
Suppose $F$ is uniformly Lipschitz with respect to $(\pi, \ell)$, continuous with respect to $\alpha$,
 and such that 
${\rm ess} \sup_{\alpha \in A}  | F( \cdot , t, 0,0,0, \alpha)|  \in \H^{p,T}$, for each $p \geq 2$. Suppose also that
 \begin{equation}\label{rrr}
F(t,\pi,l_1, \alpha)- F(t,\pi,l_2, \alpha) \geq 
\langle \tau_t^{\pi,l_1,l_2,\alpha}  \,,\,l_1 - l_2 \rangle_\nu,
\end{equation}
for some adapted process 
$\tau_t^{\pi,l_1,l_2, \alpha}(\cdot)$ 
satisfying
$
|\tau_t^{\pi,l_1,l_2, \alpha}(u) | \leq \bar \psi (u),
$
where $\bar \psi$ is bounded and in $L^p_{\nu}$, for all $p \geq 1$, and 
$\tau_t^{\pi,l_1,l_2, \alpha} \geq  -1-C_1$.

For each $\alpha \in {\cal A}$, the associated driver is given by 
\begin{equation}\label{baralpha}
 F(t, \omega, \pi, \ell, \alpha_t(\omega)).
 \end{equation}
Note that these drivers are equi-Lipschitz.
For each $\alpha$ $\in$ ${\cal A}$, let $\rho^{\alpha}$ be the dynamic risk-measure induced by the BSDE associated with $F(., \alpha_t)$ and driven by $W^{\alpha}$ and $\tilde N^\alpha$.\\
 More precisely, for each $\tau \in \T_0$ and $\zeta \in L^p({\cal F}_\tau)$ with $p>2$,
there exists a unique solution $(X^\alpha, \pi^\alpha, l^\alpha)$ in $ {\cal S}_{\alpha}^{2}\times \H_{\alpha}^{2} \times \H_{\alpha, \nu}^{2}$ of the $Q^{\alpha}$-BSDE 
 \begin{equation}\label{bsde1}
-dX^\alpha_t =F(t,\pi^\alpha_t, l^\alpha_t, \alpha_t)dt - \pi^\alpha_t dW^\alpha_t -  \int_{\R^*} l^\alpha_t(u) \tilde N^\alpha(dt,du)  ; \qquad  X^\alpha_{\tau} = \zeta,
 \end{equation}
driven by $W^{\alpha}$ and $\tilde N^\alpha$. 
 The dynamic risk-measure $\rho^{\alpha} (\zeta, \tau)$ of position $\zeta$ is thus well defined by
  \begin{equation}\label{labelroa}
  \rho^{\alpha}_t (\zeta, \tau) :=  -X^{\alpha}_t(\zeta,\tau), \;\;0\leq t \leq \tau,
  \end{equation}
 with $X^{\alpha}(\zeta,\tau)= X^{\alpha}$. Assumption (\ref{rrr}) yields the monotonicity property of $\rho^{\alpha}$.

 

The agent is supposed to be averse to ambiguity. Her dynamic risk measure is given,
for each $\tau \in \T_S$ and $\zeta \in L^p({\cal F}_{\tau})$, $p >2$,
 by 
 \begin{equation}\label{ambi}
{\rm ess} \sup_{\alpha \in {\cal A}} \rho^{\alpha}_S (\zeta, \tau) = - {\rm ess} \inf _{\alpha \in {\cal A}} X_S^{\alpha}(\zeta, \tau).
\end{equation}
at each stopping time $S$ $\in$ $\T_0$.




The financial dynamic position is given here by a RCLL predictable process $(\xi_t)$ which belongs to ${\cal S}^p$. 
At  fixed time $S \in \T_0$, the agent wants to choose a stopping time in $\T_S$ so that it minimizes (\ref{ambi}), which leads to the following mixed control/optimal stopping problem:
 \begin{equation*}
u(S):={\rm ess} \inf_{\tau \in \T_S}{\rm ess}\sup_{\alpha \in {\cal A}} \rho^{\alpha}_S (\xi_{\tau}, \tau)=- {\rm ess} \sup_{\tau \in \T_S} {\rm ess}\inf_{\alpha \in {\cal A}}  X^{\alpha}_S(\xi_{\tau}, \tau),
\end{equation*}
which corresponds to that studied in Section \ref{mixed2}.

\begin{Theorem}\label{exemple}
Let $(Y,Z,k)$ be the solution of the RBSDE associated with obstacle $(\xi_t)$ and Lispchitz driver $f$, defined for each $(t, \omega, \pi, \ell)$ by
\begin{equation}\label{d}
{f}(t, \omega, \pi, \ell):=\inf_{\alpha \in A} \{F(t, \omega, \pi, \ell, \alpha) + \beta^1 (t, \omega, \alpha) \pi + \langle \beta^2 (t, \omega, \alpha), \ell\rangle_\nu \}.
\end{equation}
 For each $S$ $\in$ $\T_0$, we have
$$Y_S=  \underline V(S)= \bar V(S)\;\; {\rm a.s.}$$ 
\end{Theorem}

\dproof In order to prove this result, we will express the problem in terms of BSDEs and RBSDEs under probability $P$ and then apply Theorem \ref{optieps}.

Fix now $\tau \in \T_0$ and $\zeta \in L^p({\cal F}_{\tau})$ with $p>2$.
Since $(X^\alpha, \pi^\alpha, l^\alpha)$ is the solution of 
BSDE (\ref{bsde1}), it clearly satisfies the following $P$-BSDE driven by $W$ and $\tilde N$
\begin{equation}\label{bsde2}
-dX^\alpha_t ={f}^\alpha(t,\pi^\alpha_t, l^\alpha_t)dt - \pi^\alpha_t dW_t -  \int_{\R^*} l^\alpha_t(u) \tilde N(dt,du)  ; \qquad  X^\alpha_{\tau} = \zeta,
 \end{equation}
where the driver is given by 
 \begin{equation}\label{barrebis}
  f^{\alpha}(t,  \pi, \ell):=F(t,  \pi, \ell, \alpha_t) +\beta^1 (t, \alpha_t) \pi +\langle \beta^2 (t, \alpha_t), \ell \rangle_\nu.
 \end{equation}
The process $(X^\alpha, \pi^\alpha, l^\alpha)$ is the solution of $P$-BSDE (\ref{bsde2}) in $ {\cal S}^2\times \H^{2} \times \H_{ \nu}^{2}$ (see the proof of Theorem 5.9 in \cite{QuenSul}). 
Moreover, for each $\alpha$, $f^\alpha$  satisfies Assumption \ref{Royer}, and $f$, defined by \eqref{d}, is a Lipschitz driver (see \cite{QuenSul}).\\
By the definition of $f$ (see \eqref{d}) and $f^\alpha$ (see \eqref{barrebis}), we get that for each $\alpha$ $\in$ ${\cal A}$, $f \leq f^{\alpha}$.

Also, 
for each $\eta >0$ and each $(t, \omega,  \pi, l)$ $\in$ $\Omega \times [0,T] \times \R \times L^2_\nu$, there exists $\alpha^{\eta}$ $\in$ $A$ such that 
$$f (t, \omega, \pi, \ell) +\eta \geq F (t, \omega,\pi, \ell, \alpha^{\eta})+ \beta^1 (t, \omega, \alpha^{\eta}) \pi + \langle 
 \beta^2 (t, \omega, \alpha^{\eta}), \ell \rangle_\nu . $$
By the section theorem of \cite{DM1}, for each $\eta >0$, there exists an $A$-valued predictable process $( \alpha^{\eta}_t)$ such that 
$ f (t,  Z_t, k_t)+\eta \geq   f ^{\alpha^{\eta}}(t, Z_t, k_t)$, $dP \otimes dt$-a.s.\,
Consequently, by  Theorem \ref{optieps}, the result follows.
\fproof

\begin{Corollary}\label{exemple2}

Suppose $A$ is compact and $F$, $\beta^1$ and $\beta^2$ are continuous with respect to $\alpha$. Suppose that the position $(\xi_t)$ is left-usc along stopping times. Then,  there exists $\bar \alpha$ $\in$ $A$ such that 
\begin{equation}\label{exist2}
f(t,Y_t, Z_t,k_t) = {\rm ess} \inf_{\alpha \in \Ac} f^{\alpha}(t,Y_t, Z_t,k_t) = 
f^{\bar \alpha}(t,Y_t, Z_t,k_t),  0 \leq t \leq T, \;\;   dt \otimes dP-a.s.
\end{equation}
Also, for each $S$ $\in$ $\T_0$, the pair $( \tau_S^*, \bar \alpha)$ is an $S$-saddle point, 
where $\tau^*_S = \inf \{u \geq S ; \, Y_u = \xi_u \}$.

This result still holds in the case when $A$, instead of being compact, is a bounded, convex and closed subset of a separable Hilbert space, and if $F$,  $\beta^1$ and $\beta^2$ are convex and lower semicontinuous with respect to $\alpha$.
\end{Corollary}
\dproof
 Since $A$ is compact and that $F$, $\beta^1$ and $\beta^2$ are continuous with respect to $\alpha$, the section theorem of \cite{DM1} provides the existence of $\bar \alpha$ $\in$ $A$ such that 
(\ref{exist2}) is satisfied. By Corollary \ref{existencedeux}, $( \tau_S^*, \bar \alpha)$ is thus an $S$-saddle point.

Let us now consider the second case. By convex analysis arguments, one can show the existence of 
$\bar \alpha$ $\in$ $A$ satisfying equality (\ref{exist2}) (for details, see the proof of Theorem 5.2 in 
\cite{QuenSul}). The result follows.
\fproof

\paragraph{Example.}
Suppose that $L^2_{\nu}$ is separable and that $A$ is a Borelian of the Hilbert space $\R \times L^2_{\nu}$ such that $A \subset [-K,K] \times \Upsilon$, where
$$
\Upsilon:= \{ \varphi \in {\cal P}, \,\, \,\,C'_1 \leq \varphi (u) \,\,\, {\rm and}\, \, |\varphi (u)| \leq \psi (u) \,\, \,\,\nu(du) \text{ a.s. }\},
$$
with $C'_1 > -1$ and $\psi $ is bounded and in $L^p_{\nu}$, for all $p \geq 1$. 
For each process $\alpha:= (\alpha^1, \alpha^2)$ $\in$ ${\cal A}$, the prior $Q^{\alpha}$ is defined as the probability measure which admits 
 $Z_{T}^{\alpha}$ as density with respect to $P$, $Z^{\alpha}$ being the solution of  
 \begin{equation*}
  dZ_t^{\alpha} = Z_{t^-} ^{\alpha}\left( \alpha^1 _t dW_t + \int_{\R^*}
\alpha^2_t(u) d\tilde N(dt,du)\right)\,; \,\,\,\, Z^{\alpha}_0= 1.
\end{equation*}
Theorem \ref{exemple} and Corollary \ref{exemple2} then hold. 
\begin{remark}
 In the case when $F(t, \omega, \pi, \ell, \alpha_t(\omega))$ is {\em linear} with respect to $\pi$ and $\ell$, the above problem is related to that studied in \cite{Bayraktar-2} (in the Brownian case). 

 \end{remark}
%

\appendix
\section{Appendix}

\begin{Proposition}\label{est}
Let $T>0$ and let $\xi$ $\in {\cal S}^2$. Let $f^1$ be a Lipschitz driver with Lipschitz constant $C$ and let  $f^2$ be a driver. For $i=1,2$, let $(Y^i, Z^i ,k^i, A^i)$  be  a solution of the RBSDE associated to 
terminal time $T$, driver $f^i$  and obstacle $\xi$.
 For $s$ in $[0,T]$, denote $\bar Y_s := Y^1_s - Y^2_s, \,\,\, \bar Z_s := Z^1_s - Z^2_s$,  $\bar k_s := k^1_s - k^2_s $, and 
 $\bar f(s): = f^1(s, Y^2_s, Z^2_s, k_s^2) - f^2(s, Y^2_s, Z^2_s, k_s^2)$.\\
Let $ \eta, \beta >0 $ be such that 
 $\beta \geq \frac{3}{\eta} +2C .$ 
If $\eta \leq \frac{1}{C^2}$, then, for each $t \in [0,T]$, we have
\begin{equation}\label{A26}
e^{\beta  t}  \bar Y_t   ^2 \leq   \eta \,E[ \int_t^T e^{\beta  s} \bar f(s) ^2  ds \mid 
{\cal F}_t ] \;\; \text{ \rm a .s.}\,\,\,\,\,{\rm and}
\end{equation}
\begin{eqnarray}\label{AA27}
\|\bar Y \|_\beta^2 \leq T  \eta
\|\bar f \|_\beta^2.
\end{eqnarray}
Also, if $\eta < \frac{1}{C^2}$, we then have 
\begin{equation}\label{AA28}
\|\bar Z \|_\beta^2 + \|\bar k \|_{\nu,\beta}^2
\leq \frac{\eta}{1 - \eta C^2}  \|\bar f \|_\beta^2.
\end{equation}
\end{Proposition}

\dproof 
From It\^o's formula applied  to the semimartingale $e^{\beta s} \bar Y_s$ 
between $t$ and $T$, it follows that
\begin{align}
e^{\beta t} \bar Y_t ^2 + \beta \int_t^T e^{\beta s} \bar Y_s^2 ds + \int_t^T e^{\beta s} \bar Z_s^2 ds &+ \int_t^T e^{\beta s} \| \bar k_s\|_\nu^2 ds \nonumber \\  
 &=  2 \int_t^T e^{\beta s} \bar Y_s (f^1 (s, Y^1_s, Z^1_s, k^1_s) - f^2 (s, Y^2_s, Z^2_s, k^2_s)) ds  \nonumber \\
 &\quad - 2 \int_t^T e^{\beta s} \bar Y_s \bar Z_s dW_s - 2 \int_t^T e^{\beta s} \int_\RB^* \bar Y_{s^-} \bar k_s(u) d \tN (du,dt) \nonumber \\
&\quad + 2 \int^T_t e^{\beta s} \overline{Y}_{s^-}  dA^1_s - 2\int^T_t e^{\beta s} \overline{Y}_{s^-}  dA^2_s \label{russ}
\end{align}
Now, we have a.s.
$$\overline{Y}_s dA^{1,c}_s = (Y^1_s - \xi_s)dA^{1,c}_s - (Y^2_s - \xi_s)dA^{1,c}_s = - (Y^2_s - \xi_s)dA^{1,c}_s\leq 0 $$
and by symmetry, $\overline{Y}_s dA^{2,c}_s  \geq 0$ a.s.\,
Also, we have a.s.
$$\overline{Y}_{s^-} \Delta A_{s}^{1,d} = (Y^1_{s^-} - \xi_{s^-}) \Delta A_{s}^{1,d} - (Y^2_{s^-} - \xi_{s^-}) \Delta A_{s}^{1,d} = - (Y^2_{s^-} - \xi_{s^-}) \Delta A_{s}^{1,d} \leq 0 $$
and $\overline{Y}_{s^-} \Delta A_{s}^{2,d}  \geq 0$ a.s.
Consequently, the two last terms of the r.h.s. of \eqref{russ} are non positive.
Moreover,
\begin{align*}
|f^1(s,Y^1_s, Z^1_s, k^1_s) - f^2(s,Y^2_s,Z^2_s,k^2_s)| & \leq |f^1(s,Y^1_s, Z^1_s, k^1_s) - f^1(s,Y^2_s,Z^2_s,k^2_s)|  + |\bar f_s| \\
& \leq C  | \bar Y_s| + (C|\bar Z_s| + C\| \bar k_s\|_\nu + |\bar f_s|).
\end{align*}
Now,  for all real numbers $y$, 
$z$, $k$, $f$ and $\varepsilon >0$\\
$ 2y (Cz + Ck  + f) \leq \frac{ y^2}{\varepsilon^2}+ \varepsilon^2(Cz+ Ck  + f)^2 \leq \frac{ y^2}{\varepsilon^2} +  3 \varepsilon^2(C^2 y^2+ C^2 k^2 +f^2)$. Hence, 
we get
\begin{align}\label{eq2a}
e^ {\beta t} & \bar Y_t^2  +  E \left[\beta \int_t^T e^{\beta s}  \bar Y_s^2 ds + \int_t^T e^{\beta s} ( \bar Z_s^2 + \| \bar k_s\|_\nu^2) ds \mid \FC_t \right] \nonumber \\
& \leq \EB \left[ (2C+\frac{ 1}{\varepsilon^2}) \int_t^T  e^{ \beta s}  \bar Y_s^2 ds + 3C^2 \varepsilon^2 \int_t^T e^{\beta s} ( \bar Z_s^2 + \| \bar k_s\|_\nu^2)ds \mid \FC_t \right] \nonumber \\
& \quad + 3 \varepsilon^2  \EB \left[ \int_t^T e^{\beta s}  \bar f_s^2 ds \mid \FC_t \right].
\end{align}
Let us make the change of variable $\eta = 3 \epsilon^2$. Then, for each  $\beta,  \eta>0$  chosen as in the theorem,  these inequalities lead to
\eqref{A26}. 
We obtain the first inequality of \eqref{AA27} by integrating  \eqref{A26}. Then \eqref{AA28}  follows from  inequality \eqref{eq2a}.
\fproof

\begin{remark} \label{AA29}
 By classical results on the norms of semimartingales, one similarly shows that
$
\|  \bar Y \|_{S^2} \leq K \|\bar f \|_{\H^2}
$, 
where $K$ is a positive constant only depending on $T$ and $C$.
\end{remark}


 
{\bf Proof of Theorem \ref{exiuni}}: Using the previous a priori estimates, we show that the mapping $\Phi$ is a contraction from ${\cal H}_\beta^2$ into ${\cal H}_\beta^2$.
Given $(U, V,l) \in {\cal H}_\beta^2$, let $ (Y, Z,k) := \Phi (U, V, l)$, that is,
the solution of the RBSDE  associated with driver process $ 
f^1_s:=f(s, U_s , V_s, l_s)$ (which does not depend on the solution).  
Let $(U', V', k' )$ be another element of ${\cal H}_\beta^2$  and let $(Y', Z',k') := 
\Phi (U', V', l')$, that is,
the solution of the RBSDE  associated with driver process $f^2_s:=
f(s, U'_s , V'_s, l'_s)$. 

Set $\bar{U} = U - U'$, $\bar{V} = V - V'$, $\bar{l} = l - l'$, $\bar{Y} = Y - Y'$, $\bar Z = Z-Z'$  $\bar{k} = k - k'$. Let $\Delta f_\cdot :=  f(\cdot, U, V, l) - f(\cdot, U', V', l')$.
Using estimates \eqref{AA27} and \eqref{AA28} with $\eta \leq \frac{1}{2C^2}$ and Lipschitz constant equal to  $0$ (since the driver $f^1$ does not depend on the solution), we get 
$$
\|\bar{Y}\|_\beta^2 + \|\bar{Z}\|_\beta^2 + \|\bar{k}\|_{\nu, \beta}^2  \leq 
\eta (T+2) \|\Delta f \|_\beta^2 \leq
\eta (T+2) 2C^2  \| \bar{U}\|_\beta^2 + \| \bar{V}\|_\beta^2 + \| \bar{l}\|_{\nu, \beta}^2),
$$
where the second inequality follows from the Lipschitz property of $f$ with constant $C$. 
Choosing 
$\eta = \frac{1} {(T+2) 4C^2 }$,
we deduce
$ \| (\overline{Y}, \overline{Z}, \overline{k}) \|_\beta^2 \leq \frac{1}{2} 
\| (\overline{U}, \overline{V}, \overline{l}) \|_\beta^2.  $
Hence, $\Phi$ is a contraction and thus admits a unique fixed point $(Y, Z, k)$ in ${\cal H}_\beta^2$, which 
corresponds to the solution of RBSDE~\eqref{RBSDE}.\\

{\footnotesize

\end{document}